\documentclass[preprint,12pt,3p]{elsarticle}
\usepackage{amsmath}
\usepackage{amssymb}
\usepackage[dvipsnames]{xcolor}
\usepackage[utf8]{inputenc}
\usepackage{verbatim}
\usepackage[mathcal]{euscript}
\usepackage{graphicx}
\usepackage[english]{babel}
\usepackage{multicol}
\usepackage[utf8]{inputenc}  
\usepackage[T1]{fontenc}  
\usepackage{caption}
\usepackage[colorlinks=true ]{hyperref}
\usepackage[final]{pdfpages}
\usepackage{multicol}
\usepackage{float}
\usepackage{tikz}
\usepackage{tikz-3dplot} 
\usepackage{pdfpages}       
\usepackage{varioref} 
\usepackage{float}
 \usepackage{multicol}
\usepackage{multirow}
\usetikzlibrary{quotes,arrows.meta}
\usetikzlibrary{shapes,arrows}
\numberwithin{equation}{section}  

\definecolor{camel}{rgb}{0.76, 0.6, 0.42}

\newtheorem{definition}{Definition}[section]

\newtheorem{Claim}[definition]{Claim}
\newtheorem{Remark}[definition]{Remark}
\newtheorem{example}[definition]{Example}

\newtheorem{proof}[definition]{Proof}
\newtheorem{Corollary}[definition]{Corollary}

\newcommand \bei {\begin{itemize}}
\newcommand \eei {\end{itemize}}
\newcommand \ubar u

\newcommand \del \partial

\newcommand \la \langle
\newcommand \ra \rangle 

\newcommand \auth    \textsc
\newcommand \be {\begin{equation}}
\newcommand \ee {\end{equation}}

\newcommand \bcor {\begin{Corollary}}
\newcommand \ecor {\end{Corollary}}

\newcommand \bpro {\begin{proof}}
\newcommand \epro {\end{proof}}
\newcommand \bdf {\begin{definition}}
\newcommand \edf {\end{definition}}
\newcommand \bex {\begin{example}}
\newcommand \eex {\end{example}}
\newcommand \bcl {\begin{Claim}}
\newcommand \ecl {\end{Claim}}
\newcommand \brm {\begin{Remark}}
\newcommand \erm {\end{Remark}}

\DeclareMathOperator{\e}{e}

\let\oldmarginpar\marginpar
\renewcommand\marginpar[1]{\-\oldmarginpar[\raggedleft\footnotesize #1]%
{\raggedright\footnotesize #1}}
\begin{document}
\begin{frontmatter}

\title{LRBF meshless methods for predicting soil moisture distribution in root zone }     

\author[label1]{Mohamed Boujoudar}
\address[label1]{ Mohammed VI Polytechnic University, Green city, Morocco}
\address[label2]{ University of Ottawa, Ottawa, Canada}
\author[label1,label2]{Abdelaziz Beljadid}

\ead{abdelaziz.beljadid@um6p.ma} 
\address[label3]{ University of Hassan II -Casablanca, Morocco}
\author[label3]{Ahmed Taik}
\begin{abstract}  
In this paper, we first propose a coupled numerical model of unsaturated flow in soils and plant root water uptake. The Richards equation and different formulations are used in the developed numerical model to describe infiltration in root zone and to investigate the impact of the plant root on the distribution of soil moisture.
The Kirchhoff transformed Richards equation is used and the Gardner model is considered for capillary pressure.
In our approach, we employ a meshless method based on localized radial basis functions (LRBF) to solve the resulting system of equations. The LRBF approach is an accurate and computationally efficient method that does not require mesh generation and is flexible in addressing high-dimensional problems with complex geometries. Furthermore, this method leads to a sparse matrix system, which avoids ill-conditioning issues.
We implement the coupled numerical model of infiltration and plant root water uptake for one, two, and three-dimensional soils. Numerical experiments are performed using nontrivial analytical solutions and available experimental data to validate the coupled numerical model. The numerical results demonstrate the performance and ability of the proposed numerical method to predict soil moisture dynamics in root zone.

\end{abstract} 
\begin{keyword}
Richards equation, Unsaturated water flow, Root water uptake, Kirchhoff transformation, Soil-water-plant interactions, Meshfree methods, Localized Radial Basis Function             
\end{keyword}
\end{frontmatter}

\section{Introduction} 
Studying water flow dynamics  is crucial in modeling hydrological processes. Water uptake by plant roots is a significant process in subsurface unsaturated flow modeling and has a significant impact on the evolution of soil moisture and nutrient transport \cite{vsimuunek2009modeling}. Predicting moisture fluxes in unsaturated soils and their interaction with vegetation has practical implications in agriculture, water management, and climate science \cite{feddes2004parameterizing}. Understanding root water uptake process is needed to estimate crops water requirements and design best water management practices. For instance, optimum irrigation scheduling for efficient water management of crops requires a proper understanding of soil, water and plant interactions \cite{green2006root,kumar2012methods}. Plant root water absorption is a complex process that involves various interactions between soil, plant, and climate \cite{vsimuunek2009modeling}. The interactions of these complex processes can be explained using robust and efficient coupled numerical models.

The Richards equation \cite{richards1931capillary} is the most commonly used model to describe water flow in unsaturated soils due to the effects of gravity and capillarity. This equation is highly nonlinear due to the nonlinear constitutive relations which describe the relationship among capillary pressure, relative permeability, and saturation \cite{gardner1958some,brooks1964hydrau,van1980closed}.

Several models have been developed to describe plant root water uptake. Molz \cite{molz1981models} classified these models into two categories. The first is based on a microscopic approach, whereas the second is based on a macroscopic one. The microscopic approach describes water extraction at an individual root level and considers soil water flow toward the particular root, predominantly radial \cite{gardner1960dynamic,molz1968soil,hillel1975microscopic}. The microscopic models are physically-based and are efficiently used in various studies \cite{molz1971interaction,hillel1975microscopic,hainsworth1986water}.

The macroscopic approach \cite{ogata1960transpiration,molz1970extraction,FeddesR.A1978Sofw} considers properties of the overall root system. In this case, the root water extraction is considered as a volumetric sink term in the Richards equation. The macroscopic models are typically empirical \cite{FeddesR.A1978Sofw,genuchten1987numerical,albasha2015compensatory}. Several macroscopic models have been developed based on different approaches. For instance, Gardner \cite{gardner1964relation} proposed a macroscopic model describing the dynamic aspects of water availability to plants. Feddes et al. \cite{FeddesR.A1978Sofw} proposed a relatively simple empirical model to describe root water uptake. Molz and Remsen \cite{molz1970extraction} proposed a linear model where the extraction rate is assumed to vary linearly with the soil depth. Prasad \cite{prasad1988linear} proposed another linear model for root water uptake which depends on the root depth and rate of
evapotranspiration.  In \cite{prasad1988linear}, the author showed that the linear models give satisfactory results in predicting soil moisture in the root zone. Li et al. \cite{li1999exponential} suggested an exponential root water uptake model derived from measured root length distributions. In comparison with field data, they showed that the exponential model \cite{li1999exponential} performs well compared to constant and linear models \cite{FeddesR.A1978Sofw,prasad1988linear}. Ojha and Rai \cite{ojha1996nonlinear} proposed a nonlinear root water extraction model and investigated its performance on irrigation scheduling compared to experimental data \cite{erie1965consumptive}. Several studies have favoured the use of macroscopic models for their simplicity since they don't require information about the detailed geometry of the root system. For more details, we refer to the previous studies  \cite{molz1970extraction,molz1981models,mathur1999modeling,de2012root}.

Incorporating the plant root water uptake as a sink term in the Richard equation adds more complexity to its analytical and numerical solutions. Except for some simple cases, analytical solutions have been developed \cite{raats1974steady, raats1976analytical,lomen1976solution, rubin1993stochastic, broadbridge1999forced, basha2000multidimensional, schoups2002analytical, yuan2005analytical}. For instance, Lomen and Warrick  \cite{lomen1976solution} presented analytical solutions for the Richards equation under steady-state conditions. Their solutions were based on the Gardner model for capillary pressure, and they incorporated a sink term that depends on the matrix flux potential. Yuan and Lu \cite{yuan2005analytical} developed exact solutions using the Gardner model \cite{gardner1958some} for the soil hydraulic conductivity and different formulations for root water uptake sink term. Broadbridge et al. \cite{broadbridge2017exact} developed exact solutions for the Richards equation with a sink term nonlinearly dependent on soil water content. In their solutions, they used the Gardner model for the capillary pressure. Most of these analytical solutions assume a simplified model for root water uptake and consider specific initial and boundary conditions. In more practical situations, efficient numerical approaches in terms of computing cost and accuracy are required.

Various numerical approaches have been proposed to solve the Richards equation without a sink term \cite{celia1990general,boujoudar2021localized,boujoudar2021modelling,keita2021implicit,boujoudar2023localized}. Nevertheless, the design of appropriate numerical schemes for the Richards equation that includes a sink term remains challenging. Few numerical techniques have been proposed to solve the Richards equation considering root water uptake as a sink term. For instance, Neuman et al. \cite{neuman1975finite} used a finite-element method to solve the soil water flow equation, taking into account water uptake by roots. Feddes et al. \cite{feddes1976simulation} proposed an implicit finite-difference method to solve the Richards equation considering plant root water uptake as sink term. Vrugt et al. \cite{vrugt2001one} introduce the Galerkin finite element method to solve the Richards equation with a sink term of water uptake by roots. Rees and Ali \cite{rees2006seasonal} used the finite-element method to solve the Richards equation incorporating a sink term for root water extraction. Difonzo et al. \cite{difonzo2021shooting} proposed a numerical approach based on the shooting method to solve the unsaturated flow equation with root water uptake models. Despite the studies conducted, there remains a need for more advanced numerical methods for computing water flow in unsaturated soils taking into account the absorption of water by the roots.

This study presents a coupled numerical model that accounts for both unsaturated soil flow and plant root water uptake. The proposed model uses the Richards equation and includes different root water uptake formulations as sink terms in the governing equation. The first formulation is proposed by Yoan et al. \cite{yuan2005analytical}, where the stepwise and exponential functional forms are considered for the root water uptake. The second formulation is proposed by Broadbridge et al. \cite{broadbridge2017exact}, where root water uptake's sink term is assumed nonlinearly dependent on soil water content. In all considered formulations, the Gardner model \cite{gardner1958some} is employed for the capillary pressure.

In our approach, we employ a meshless method based on localized radial basis function (LRBF) to solve the resulting system of equations. The RBF meshless approach is an accurate and computationally efficient method that eliminates the need for mesh generation and is flexible in addressing high-dimensional problems with complex geometries. 
RBF meshless methods can be categorized into global \cite{kansa1990multiquadrics,kansa1990multiquadrics1} and local \cite{lee2003local,li2013localized} methods. In the global approach, the collocation is performed globally over the entire computational domain, while the local approach performs collocation locally over a set of influence domains. Although the global method is simple to implement \cite{kansa1990multiquadrics,kansa1990multiquadrics1}, it faces two challenges which are: ill-conditioned resultant matrices and the issue of selecting an appropriate shape parameter for certain radial basis functions \cite{sundin2020global}. To overcome these drawbacks, LRBF methods are proposed \cite{lee2003local,vsarler2007global,li2013localized}. These local approaches result in a sparse matrix that avoids ill-conditioning issues encountered in the global approach, where a full matrix is produced \cite{li2013localized}. LRBF meshless methods have been effective in solving various partial differential equations, such as the Richards equation without the sink term \cite{stevens2009meshless,stevens2011scalable,ben2018radial,boujoudar2021localized,boujoudar2021modelling}. To our best knowledge, no studies have explored the use of LRBF methods in solving the Richards equation that includes a sink term accounting for root water uptake.

In this paper, we propose LRBF meshless method \cite{li2013localized} to solve the Richards equation with the presence of a sink term of water uptake by plant roots. To deal with the nonlinearity of the governing equation, we employed the Kirchhoff transformation. 
To linearize the nonlinear problem resulting from the Kirchhoff transformation technique, Picard's iterations are used.
The numerical results are validated and compared with nontrivial analytical solutions and available experimental data. 

The structure of the paper is organized as follows: Section \ref{sec:2} introduces the coupled mathematical model of infiltration and plant root water uptake. Section \ref{sec:3} presents the proposed LRBF meshless method to solve the governing system. Numerical experiments for modeling soil moisture distribution in the root zone are performed in Section \ref{sec:4}. Conclusions are provided in Section \ref{sec:5}.
\section{Mathematical model}
\label{sec:2}

\subsection{The governing equations} 
\label{sec:2.1}
The Richards equation~\cite{richards1931capillary}, which governs water flow in the unsaturated zone while accounting for root water uptake with a sink term, is given by:
 \begin{equation}\label{E1}
  \dfrac{\partial\theta}{\partial t}-\nabla.\left(D(\theta\right)\nabla \theta)-\dfrac{\partial K(\theta)}{\partial z}= -R(\boldsymbol{x},\theta),  \text{ $\boldsymbol{x}\in\Omega, 0< t\leq T $ },
  \end{equation}
 where $\theta$ $ [L^{3}/L^{3}]$ is the volumetric water content, $D=K d h/d\theta$ $ [L^{2}/T]$ is the soil-water diffusivity, $h$ $[L]$ is the pressure head, $K$ $ [L/T]$ is the unsaturated hydraulic conductivity which is given by $K=K_s k_r$ where $K_s$ $[L/T]$ is the saturated hydraulic conductivity and $k_r$ $ [-]$ is the water relative permeability, $R$ $[T^{-1}]$ is the root water uptake, the coordinate vector $\boldsymbol{x}=(x, y, z)^{T}$ consists of the lateral directions $x$ $[L]$ and $y$ $[L]$, and the vertical direction $z$ $[L]$, which is positive upward, $\Omega $ is an open bounded set of $\mathbb{R}^{d}$, where $d$ represents the dimension of the computational domain and $T$ is a fixed time. 
 
 The Richards equation is highly nonlinear due to the nonlinear constitutive
relations which describe the relationship among capillary pressure, relative permeability and
saturation \cite{gardner1958some,van1980closed,brooks1964hydrau}. The inclusion of the root-water uptake term in the Richards equation makes it further complicated in terms of numerical resolution. In this study, we simplify our analysis by considering a homogeneous soil. We adopt the assumption made by Broadbridge et al. \cite{broadbridge2017exact}, wherein the derivative of $K$ with respect to $\theta$ is expressed as follows:
\begin{equation}\label{E3}
\frac{d K}{d \theta}=\alpha D,
\end{equation}
where $\alpha $ $ [1/L]$ is a pore size parameter, then the Richards equation becomes:
  \begin{equation}\label{E4}
  \dfrac{\partial\theta}{\partial t}-\nabla.(D(\theta)\nabla \theta)-\alpha D\dfrac{\partial \theta}{\partial z}= -R(\boldsymbol{x},\theta).
  \end{equation}
\subsection{Root water uptake models}
The sink term $R$ represents the rate of water uptake by plants, expressed as the volume of water removed per unit of time from a unit volume of soil. Different mathematical models are developed to describe plant root water uptake, which are classified into two approaches: microscopic \cite{gardner1960dynamic,molz1968soil,hillel1975microscopic} and macroscopic \cite{hoogland1981root,feddes1982simulation,li1999exponential,molz1970extraction} approaches. Here, different macroscopic models are considered sink terms in the Richards equation to study the impact of plant root absorption on soil moisture distribution. 
\subsubsection{Model I }
\begin{figure}[ht!]
\centering
\includegraphics[width=7cm,height=6cm,angle=0]{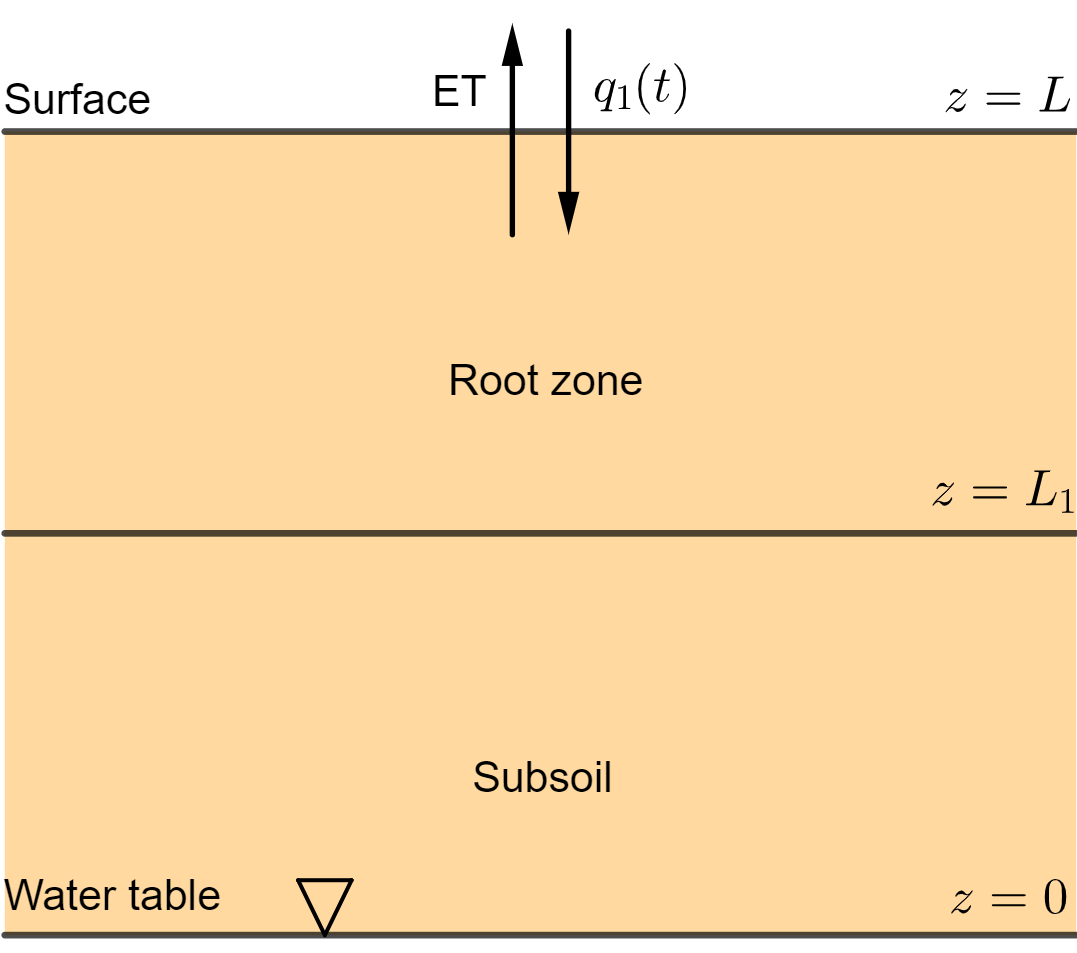}  
\caption{Schematic of vadose zone.}\label{YoanP}
\end{figure}

We present the first model used in \cite{yuan2005analytical} where the stepwise and exponential functional forms are considered for the root water uptake. In this case, the Gardner model \cite{gardner1958some} for the soil hydraulic conductivity function and a simple formulation \cite{warrick2003soil} for the volumetric water content are used:
\begin{equation}\label{E5}
K=K_s \exp(\alpha h),
\end{equation}
\begin{equation}\label{E6}
\theta=\theta_r+(\theta_s-\theta_r) \exp(\alpha h),
\end{equation}
where $\theta_{r}$ $[L/L]$ and $\theta_{s}$ $[L/L]$ are the residual and saturated moisture content respectively. We consider an unsaturated soil zone stretching from the soil surface ($z = L$) to the water table ($z = 0$). As shown in Figure \ref{YoanP}, $L_1$ denotes the maximum root depth $[L]$, $q_1$ is time-dependent flux applied at the soil surface $[L/T]$ and $ET$ is the evapotranspiration through the root zone $[L/T]$.
Assuming $R_0$ represents the maximum uptake at the land surface $[T^{-1}]$, the stepwise form of the root water uptake may be written as \cite{yuan2005analytical}:
\begin{equation}\label{E7}
    R(z)=R_0\delta(z-L_1),
\end{equation}
where $\delta(z-L_1)$ is the Heaviside function given by: 
\begin{equation}\label{E8}
\delta(z-L_1)=
\begin{cases}
    1, & \text{if}~ L_1\leq z \leq L, \\
    0, & \text{if} ~0\leq z < L_1. \\
    \end{cases}
\end{equation}
Root water extraction can be represented using the exponential form \cite{raats1974steady,rubin1993stochastic,yuan2005analytical} given by:
\begin{equation}\label{E9}
    R(z)=R_0\exp\left(\beta(z-L)\right),
\end{equation}
where $\beta$ $[L^{-1}]$ is a parameter represents the rate of reduction of root uptake. 
\subsubsection{Model II}
The second model is proposed by Broadbridge et al. \cite{broadbridge2017exact}. In this case, water uptake by plant roots is assumed nonlinearly dependent on soil water content. The soil diffusivity, unsaturated hydraulic conductivity, and root water uptake rate are expressed as follows \cite{broadbridge2017exact}:  
\begin{equation}\label{E10}
    \begin{cases}
    D=\dfrac{1}{\alpha^2 t_s}\dfrac{m}{(\e^{m}-1)}\e^{m\Theta},\\
    K=K_s\dfrac{\e^{m\Theta}-1}{(\e^{m}-1)},\\
    R=\left(\dfrac{\theta_s}{t_s}\right)\left( \dfrac{-k}{(e^{m}-1)}(e^{\Theta m}-1)-\dfrac{A}{m}(1-e^{-\Theta m})\right),
    \end{cases}
\end{equation}
where $\Theta=(\theta-\theta_r)/(\theta_{s}-\theta_r)$ $[-]$ is the normalized water content, $t_s=\theta_s/(\alpha K_s)$ $[T]$ is the gravity time scale, $A$ $[-]$, $m$ $[-]$ and $k$ $[-]$ are constant parameters depend on plant root water uptake rate and verify \cite{broadbridge2017exact}:
  \begin{equation}\label{E11}
      \vert k \vert <\dfrac{-A}{m}\e^{-m}(1-\e^{-m}),~ A<0~\text{and}~m>0.
  \end{equation}
 At high moisture contents, $R$ is very high and close to a maximum value called the potential extraction rate $R_s$ $[T^{-1}]$. As soil moisture decreases and is near the wilting point, $R$ decreases as well and becomes very low and close to zero. 
\subsection{Kirchhoff-transformed gouverning equation}
We use the Kirchhoff transformation to circumvent some difficulties encountered in solving Richards' equation due to its high nonlinearity. This transformation is given by:
\begin{equation}\label{E12}
\mu=\int_{0}^{\theta}D(\theta)d\theta.
\end{equation}
According to models presented before, this integral can be analytically calculated, which results in the following:
\begin{equation}\label{E13}
\mu=
    \begin{cases}
    \dfrac{K_s }{(\theta_s-\theta_r)}\theta, & \text{Model I},\\
    \left(\dfrac{\theta_s-\theta_r}{\alpha^2t_s}\right)
    \left(\dfrac{\e^{m\Theta}-1}{\e^{m}-1}\right), & \text{Model II}.
    \end{cases}
\end{equation}
Using this transformation, Equation \eqref{E4} can be simplified to:
\begin{equation}\label{E14}
  \dfrac{1}{D(\mu)}\dfrac{\partial\mu}{\partial t}-\nabla^{2}.\mu-\alpha \dfrac{\partial \mu}{\partial z}=-R(\mu).
\end{equation}
The boundary conditions are transformed according to the Kirchhoff variable: \\For Dirichlet boundary conditions, we assume $\theta=\theta_d$, which leads to $ \mu=\mu(\theta_d)$. In the case of Neumann boundary conditions, two cases are considered. First, when $-K\dfrac{\partial h}{\partial x^{i}}=q_x$, it follows that $-\dfrac{\partial \mu}{\partial x^{i}}=q_x$. Second, when $-K\left[\dfrac{\partial h}{\partial z}+1\right]=q_z$, it implies that $-\left[\dfrac{\partial \mu}{\partial z}+\alpha\mu\right]=q_z$, where, $\theta_d$, $q_x$ and $q_z$ are given functions associated with the boundary conditions. The variable $i$ can take the values of ${1,2}$, which correspond to the horizontal spatial dimensions $(x^{1},x^{2})=(x,y)$.
The high nonlinearity of the governing equation is reduced using the Kirchhoff transformation, which has the advantage of improving the convergence properties of the proposed numerical method \cite{boujoudar2021modelling, boujoudar2023localized, stevens2011scalable}. Once the approximate values of $\mu$ are determined by solving Equation \eqref{E14}, we can deduce the estimated values of $\theta$ based on Equation \eqref{E13}:
\begin{equation}\label{E15}
\theta=
    \begin{cases}
    \dfrac{(\theta_s-\theta_r)} {K_s}\mu, & \text{Model I},\\
    \theta_r+\dfrac{(\theta_s-\theta_r)}{m}\log\left(1+\alpha^2t_s\left(\dfrac{\e^{m}-1}{\theta_s-\theta_r}\right)\mu\right), & \text{Model II}.
    \end{cases}
\end{equation}
Evenly, one can deduce the approximate solutions of $h$ and $\Theta$ using the constitutive relations given by Equations \eqref{E6} and \eqref{E10}.
\section{Materiel and methodology }
\label{sec:3}
In this section, we present the numerical approach proposed to solve Equation \eqref{E14}. Specifically, we employ LRBF meshless techniques to solve the Richards equation, which includes a sink term to account for the root water uptake. For temporal discretization, we utilize an implicit scheme based on the backward differentiation formula (BDF). This choice ensures stability when solving Richards' equation and permits the use of reasonable time steps, which enhance the computational efficiency. To linearize the system, Picard's iterations are utilized.  
 \subsection{Linearization techniques}
 The backward Euler method of Equation \eqref{E14} is given by:
\begin{equation}\label{E16}
  \dfrac{1}{D^{n+1}}\dfrac{\mu^{n+1}-\mu^{n}}{ \Delta t}-\nabla^{2}.\mu^{n+1}-\alpha\dfrac{\partial \mu^{n+1}}{\partial z}= -R^{n+1},
\end{equation}
where $n\geq 0$, $\Delta t $ is the time step, $\mu^{n+1}$, $D^{n+1}$ and $R^{n+1}$ are the estimate values of $\mu$, $D$ and $R$ at time $t^{n+1}=(n+1)\Delta t$, respectively. Starting from a given initial condition $\mu^{0}=\mu(\theta_{\{t=0\}})$, Equation \eqref{E16} is linearized using the Picard iteration scheme \cite{celia1990general}, which is given by: 
\begin{equation}\label{E17}
  \dfrac{1}{D^{n+1,m}}\dfrac{\mu^{n+1,m+1}-\mu^{n}}{ \Delta t}-\nabla^{2}.\mu^{n+1,m+1}-\alpha \dfrac{\partial \mu^{n+1,m+1}}{\partial z}= -R^{n+1,m},
\end{equation}
where $m$ denotes iteration levels and $\mu^{n+1,0}=\mu^{n}$ is the initial guess of $\mu$. At $(m+1)^{th}$ iteration, Equation \eqref{E16} is solved for $\mu^{n+1,m+1}$ using the values of $R^{n+1,m}$ and $D^{n+1,m}$ obtained from the previous iteration. Subject to  boundary and initial conditions, Equation \eqref{E17} can be written as:
\begin{equation}\label{E18}
\begin{cases}
\mathcal{L}^{m}\mu^{n+1,m+1}(\boldsymbol{x})=f^{n,m}(\boldsymbol{x}), & \text{ $\boldsymbol{x}\in\Omega$},\\
\mathcal{B}\mu^{n+1,m+1}(\boldsymbol{x})=g(\boldsymbol{x}), & \text{ $\boldsymbol{x}\in\partial\Omega$},\\
\mu^{0}(\boldsymbol{x})=\mu_{0}(\boldsymbol{x}), & \text{ $\boldsymbol{x}\in\Omega$},
\end{cases}
\end{equation}
where $\mathcal{L}^{m}$ and $f^{n,m}$ are given by:
\begin{equation}\label{E19}
\mathcal{L}^{m}=\left( \dfrac{1}{D^{n+1,m}}\dfrac{1}{\Delta t}.-\nabla^{2}.-\alpha\dfrac{\partial . }{\partial z}\right) ,
\end{equation}
\begin{equation}\label{E20}
f^{n,m}=\left( \dfrac{\mu^{n}}{D^{n+1,m}}\dfrac{1}{\Delta t}-R^{n+1,m}\right).
\end{equation}
 $\mathcal{L}^{m}$ denotes the differential operator, which is linearized during each Picard iteration $(m+1)$, $\mathcal{B}^{m}$ is a linear border operator, which can be either Dirichlet or Neumann. At each time level $n$ and iteration $(m+1)$, LRBF meshless method \cite{lee2003local} is used to solve Equation \eqref{E18} until the stop condition is satisfied, which is determined by:
\begin{equation}\label{E21}
\vert\delta^{n+1,m+1}\vert=\vert \mu^{n+1,m+1}-\mu^{n+1,m} \vert<TOL,
\end{equation}
where $TOL$ is the error tolerance. 
\subsection{Meshless method based on LRBF
}
\begin{figure}[ht!]
\centering
\includegraphics[width=5.5cm,height=5.5cm,angle=0]{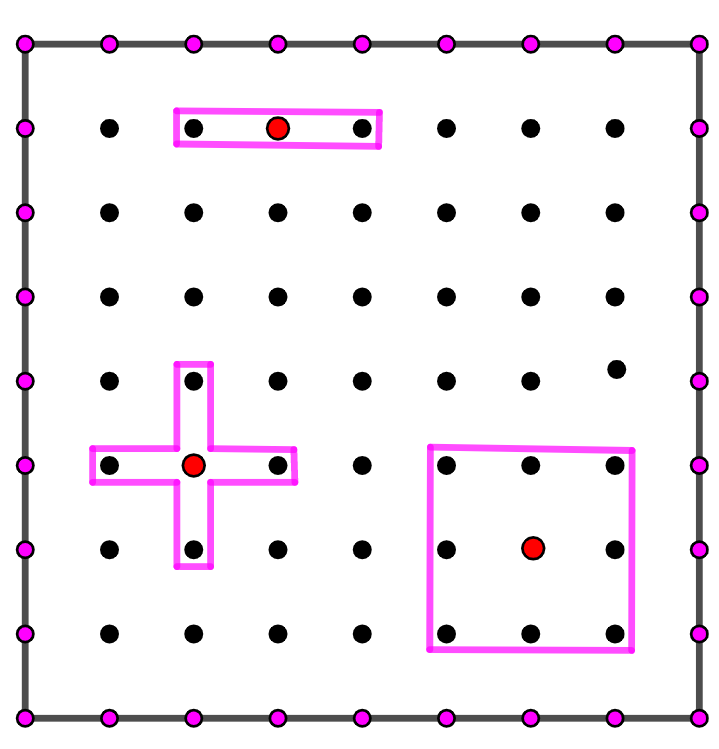}  
\caption{Schematic of influence domains $\Omega^{[s]}$ for $n_s=3,~5,~\text{and}~9$. }\label{f2}
\end{figure}

Let $\left\lbrace {\boldsymbol{x}_{j}}\right\rbrace_{j=1}^{N_i}$ be a set of $N_i$ distinct interior collocation points located in $\Omega$, and $\left\lbrace \boldsymbol{x}_{j} \right\rbrace _{j=N_i+1}^{N}$ be the boundary points, where $N$ is the total number of collocation points in $\Omega\cup\partial\Omega$. For each point $\boldsymbol{x}_{s}\in\Omega$, $s=1,2,...,N$, a localized influence domain $\Omega^{[s]}$ is created using the kd-algorithm \cite{bentley1975multidimensional}. It contains $n_{s}$ nearest neighbors interpolation points $ \left\lbrace {\boldsymbol{x}^{[s]}_{j}}\right\rbrace _{j=1}^{n_s} $ to $\boldsymbol{x}_{s}$. Figure \ref{f2} shows three examples of local influence domains including $n_s=3,~5~\text{and}~9$.
According to previous research \cite{li2013localized,hamaidi2016space,boujoudar2023localized}, the number of neighboring points $n_s$ needed for an accurate scheme may vary depending on the dimension of the computational domain. For instance, Lee et al. \cite{lee2003local} recommend to select 
 $2\dim(\Omega)$+1 nodes in each influence domain. The approximate solution of $\mu^{m+1,n+1}_{[s]}$ can be expressed as a linear combination of RBFs within each influence domain:
\begin{equation}\label{E22}
\mu^{n+1,m+1}_{[s]}(\boldsymbol{x}_{s})=\sum_{k=1}^{n_{s}}\lambda^{[s]^{n+1,m+1}}_{k}\varphi(\Vert \boldsymbol{x}_{s}-\boldsymbol{x}_{k}^{[s]} \Vert),
\end{equation}
where $\{\lambda_{k}^{[s]^{n+1,m+1}}\}_{k=1}^{n_s}$ are unknown coefficients and $\varphi$ is a RBF. In our case, we employ the exponential RBF defined as:
\begin{equation}\label{E23}
    \varphi(r_k)=\exp(-\varepsilon^2r_k^2),
\end{equation}
where $r_k=\Vert \boldsymbol{x}_{s}-\boldsymbol{x}_{k}^{[s]} \Vert$ indicates the distance between $\boldsymbol{x}_{s}$ and $\boldsymbol{x}_{k}^{[s]}$, and $\varepsilon>0$ is an arbitrary shape parameter. The exponential RBF is chosen due to its proven positive definiteness, as demonstrated in various studies \cite{micchelli1984interpolation,cheng2012multiquadric,fasshauer2007meshfree}. This characteristic is crucial to guarantee that the resulting matrix will be non-singular \cite{micchelli1984interpolation,cheng2012multiquadric,fasshauer2007meshfree,musavi1992training,garmanjani2018rbf}. 
The choice of an appropriate shape parameter is crucial for the accuracy of RBF meshless methods \cite{fasshauer2007choosing}. Various studies have proposed optimal choices for the shape parameter \cite{hardy1971multiquadric,franke1998solving}.
In this study, the use of LRBF methods allows us to overcome the issue of selecting the optimal shape parameter, as shown in previous studies \cite{lee2003local,khoshfetrat2013numerical}, they are less sensitive to the choice of shape parameter compared to the global RBF methods. From Equation \eqref{E22}, we obtain the matrix form of the solution $\mu_{[s]}^{n+1,m+1}$:
\begin{equation}\label{E24}
\mu_{[s]}^{n+1,m+1}=\varphi^{[s]}\lambda^{[s]^{n+1,m+1}},
\end{equation}
where $\mu_{[s]}^{n+1,m+1}=\left[ \mu_{[s]}^{n+1,m+1}(\boldsymbol{x_{1}^{[s]}}), \mu_{[s]}^{n+1,m+1}(\boldsymbol{x_{2}^{[s]}}),..., \mu_{[s]}^{n+1,m+1}(\boldsymbol{x_{n_{s}}^{[s]}}) \right]^{T} $,\\ $\lambda^{[s]^{n+1,m+1}}=\left[ \lambda^{[s]^{n+1,m+1}}(\boldsymbol{x_{1}^{[s]}}), \lambda^{[s]^{n+1,m+1}}(\boldsymbol{x_{2}^{[s]}}),..., \lambda^{[s]^{n+1,m+1}}(\boldsymbol{x_{n_{s}}^{[s]}}) \right]^{T} $ and $\varphi^{[s]}=\left[ \varphi(\Vert \boldsymbol{x}_{i}^{[s]}-\boldsymbol{x}_{j}^{[s]} \Vert ) \right]_{1\leq i,j\leq  n_s}$ is an $n_{s}\times n_{s}$ real symmetric coefficient matrix. The unknown coefficients $\lambda^{[s]^{n+1,m+1}}$ can be obtained as follow:
\begin{equation}\label{E25}
\lambda^{[s]^{n+1,m+1}}=(\varphi^{[s]})^{-1}\mu^{n+1,m+1}_{[s]}.
\end{equation}
Using the linear operator $\mathcal{L}^{m}$ to Equation~(\ref{E22}) yields the following equations for $\boldsymbol{x_{s}} \in \Omega$:
\begin{equation} \label{E26}
\begin{split}
\mathcal{L}^{m}\mu_{[s]}^{n+1,m+1}(\boldsymbol{x_{s}}) & = \sum_{k=1}^{n_{s}}\lambda^{[s]^{n+1,m+1}}_{k}\mathcal{L}^{m}\varphi(\Vert \boldsymbol{x_{s}}-\boldsymbol{x_{k}}^{[s]} \Vert )=\sum_{k=1}^{n_{s}}\lambda^{[s]^{n+1,m+1}}_{k}\Psi^{m}(\Vert \boldsymbol{x_{s}}-\boldsymbol{x_{k}}^{[s]} \Vert ) \\ 
& = \vartheta_{[s]}^{m}\lambda^{n+1,m+1}_{[s]}=\vartheta_{[s]}^{m}(\varphi^{[s]})^{-1}\mu^{n+1,m+1}_{[s]}=\Lambda_{[s]}^{m}\mu_{[s]}^{n+1,m+1},
\end{split}
\end{equation}
where $\Psi^{m}=\mathcal{L}^{m}\varphi$, $\vartheta_{[s]}^{m}=\left[ \Psi(\Vert \boldsymbol{x_{s}}-\boldsymbol{x_{1}}^{[s]} \Vert),...,\Psi(\Vert \boldsymbol{x_{s}}-\boldsymbol{x_{n_{s}}^{[s]}} \Vert) \right]$ and $ \Lambda_{[s]}^{m}=\vartheta_{[s]}^{m}(\varphi^{[s]})^{-1} $. In order to reformulate Equation \eqref{E26} in terms of the global $\mu^{n+1,m+1}$ instead of the local $\mu_{[s]}^{n+1,m+1}$, $\Lambda^{m}$ is considered as the expansion of $\Lambda_{[s]}^{m}$ obtained by padding the local vector with zeros at the proper positions. It follows that:
\begin{equation} \label{E27}
    \mathcal{L}^{m}\mu_{[s]}^{n+1,m+1}(\boldsymbol{x_{s}}) =\Lambda^{m} \mu^{n+1,m+1},
\end{equation}
where $ \mu^{n+1,m+1}=\left[ \mu^{n+1,m+1}(\boldsymbol{x_{1}}),\mu^{n+1,m+1}(\boldsymbol{x_{2}}),...,\mu^{n+1,m+1}(\boldsymbol{x_{N}}) \right]^{T}$. In the same way, we apply the linear operator $\mathcal{B}$ to Equation \eqref{E22} for $ \boldsymbol{x_{s}}\in\partial\Omega$:
\begin{equation} \label{E28}
\begin{split}
\mathcal{B}\mu_{[s]}^{m+1,n+1}(\boldsymbol{x_{s}})&=\sum_{k=1}^{n_{s}}\lambda^{[s]^{m+1,n+1}}_{k}\mathcal{B}\varphi(\Vert \boldsymbol{x_{s}}-\boldsymbol{x_{k}}^{[s]} \Vert)=(\mathcal{B}\varphi^{[s]})\lambda^{[s]^{m+1,n+1}}\\
&=(\mathcal{B}\varphi^{[s]})(\varphi^{[s]})^{-1}\mu_{[s]}^{n+1,m+1}=\mathbb{\sigma}^{[s]}\mu_{[s]}^{n+1,m+1}=\mathbb{\sigma}\mu^{n+1,m+1},
\end{split}
\end{equation}
where $\mathbb{\sigma}^{[s]}$ is defined as $(\mathcal{B}\varphi^{[s]})(\varphi^{[s]})^{-1}$, and $\mathbb{\sigma}$ is the expansion of $\mathbb{\sigma}^{[s]}$ that is obtained by
adding zeros at the appropriate positions. Combining Equations \eqref{E27} and \eqref{E28} into Equation \eqref{E18}, we obtain the following system:
\begin{equation}\label{E29}
\begin{split}
\mathcal{L}^{m}\mu^{n+1,m+1}(\boldsymbol{x_{s}})=\Lambda^{m}(\boldsymbol{x_{s}})\mu^{n+1,m+1}=f^{m,n}(\boldsymbol{x_{s}}),\\
\mathcal{B}\mu^{n+1,m+1}(\boldsymbol{x_{s}})=\mathbb{\sigma}(\boldsymbol{x_{s}})\mu^{n+1,m+1}=g(\boldsymbol{x_{s}}).
\end{split}
\end{equation}
As a result, we obtain a sparse system: 
\begin{equation}\label{E30}
\left(
\begin{array}{c}
\Lambda^{m}(\boldsymbol{x_{1}}) \\
\Lambda^{m}(\boldsymbol{x_{2}}) \\
 . \\
.  \\
\Lambda^{m}(\boldsymbol{x_{N_{i}}})\\
\mathbb{\sigma}(\boldsymbol{x_{N_{i}+1}}) \\
 . \\
.  \\
\mathbb{\sigma}(\boldsymbol{x_{N}})
\end{array}
\right)\left(
\begin{array}{c}
\mu^{n+1,m+1}(\boldsymbol{x_{1}})\\
\mu^{n+1,m+1}(\boldsymbol{x_{2}})\\
.\\
.\\
\mu^{n+1,m+1}(\boldsymbol{x_{N_{i}}})\\
\mu^{n+1,m+1}(\boldsymbol{x_{N_{i}+1}})\\
. \\
.\\
\mu^{n+1,m+1}(\boldsymbol{x_{N}})
\end{array}
\right)=\left(
\begin{array}{c}
f^{n,m}(\boldsymbol{x_{1}})\\
f^{n,m}(\boldsymbol{x_{2}})\\
.\\
.\\
f^{n,m}(\boldsymbol{x_{N_{i}}})\\
g(\boldsymbol{x_{N_{i}+1}})\\
. \\
.\\
g(\boldsymbol{x_{N}})
\end{array}
\right).
\end{equation}
This localized approach leads to inverting a sparse matrix, which avoids ill-conditioning problems that occur in the full matrix generated using the global method. By solving this sparse system, we get the approximate solution $\mu^{n+1,m+1}$ at all given points. Once the condition given by Equation \eqref{E21} is satisfied, we assign $\mu^{n+1}=\mu^{n+1,m+1}$.
\section{\textbf{Numerical simulations}}\label{sec:4}

This section presents numerical experiments to solve the Richards equation with a sink term, which represents the plant root water uptake. Our approach is based on LRBF meshless method based on the exponential function to solve the governing system \eqref{E20}. The first numerical test investigates soil water content dynamics during evaporation process. In the second example, root water uptake is described using model I. The third numerical test focuses on a 2D irrigation furrows system and aims to explore the influence of plant roots uptake on soil moisture distribution. The final numerical test applies the proposed numerical model to a 3D irrigation system from a circular source, evaluating its ability to predict 3D soil moisture distribution in the root zone. The numerical results are validated in comparison with non trivial analytical solutions and available experimental data. Note that, in the numerical tests where $z$ is considered positive downward, the same equation \eqref{E1} is considered by replacing $z$ with $-z$.

\subsection{Test 1: Soil water dynamics during evaporation process}
\begin{figure}[ht!]
\centering
\begin{tabular}{cc}\includegraphics[width=9cm,height=5.5cm,angle=0]{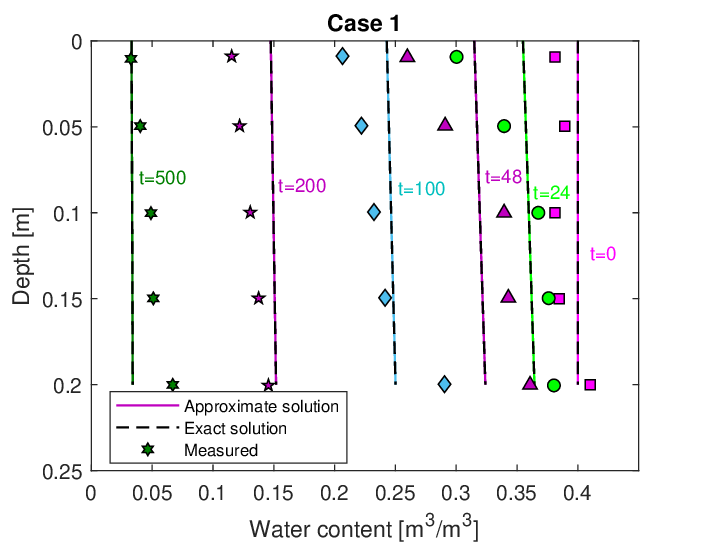}&
\includegraphics[width=7.5cm,height=5.5cm,angle=0]{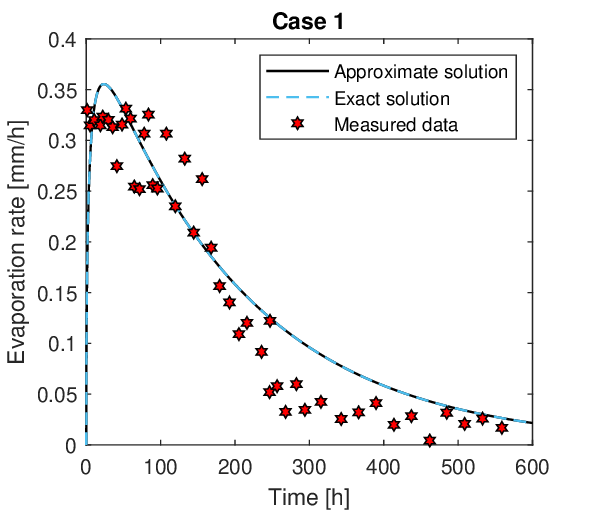}  \\
\includegraphics[width=9cm,height=5.5cm,angle=0]{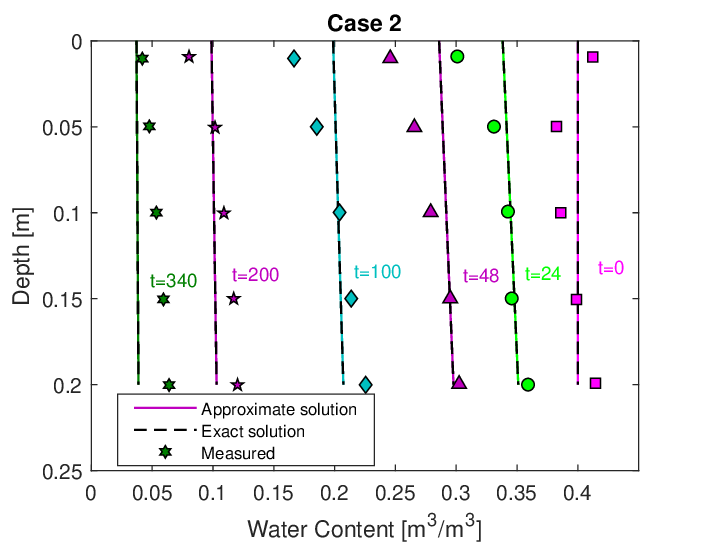} &\includegraphics[width=7.5cm,height=5.5cm,angle=0]{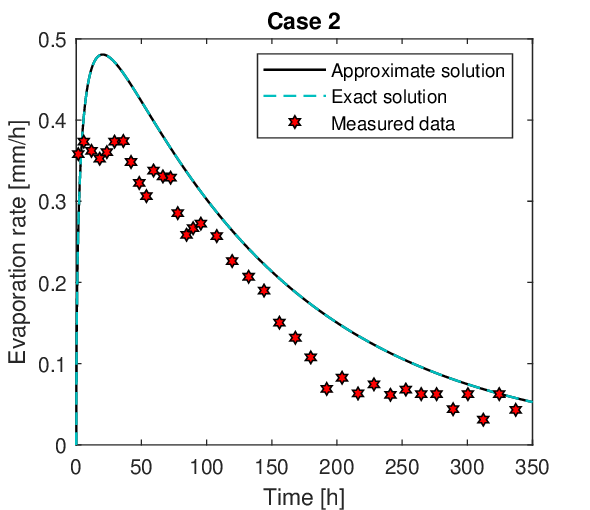} 
\\
\includegraphics[width=9cm,height=5.3cm,angle=0]{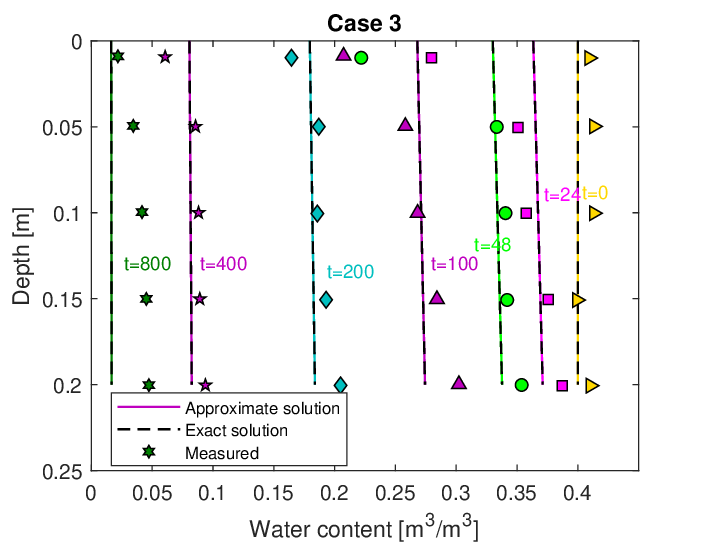}  &\includegraphics[width=7.5cm,height=5.3cm,angle=0]{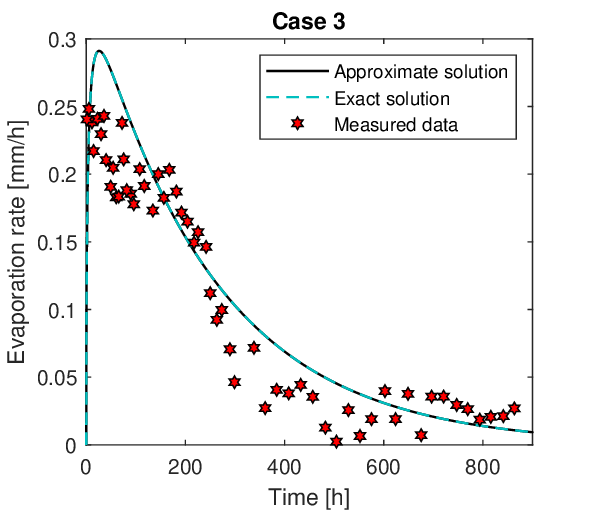}
\end{tabular}
\caption{Left: measured and computed water content profile for three cases. Right: measured and computed evaporative rate. Symbols present the experimental profile, dashed lines are for exact solution and the solid lines are for approximate solution.}\label{evapTest}
\end{figure}

In this numerical test, we consider a laboratory evaporation experience which is carried out by Teng et al. \cite{teng2012climate,teng2013analytical}. A kind of soil in Japan named K-7 sand was selected in this experience to investigate soil water content dynamics during evaporation. The soil was filled into a cylinder with a height of $20$ cm and a diameter of $10$ cm and then was wetted to saturation. Three soil samples were subjected to different environmental conditions which include relative humidity, wind speed and temperature. These conditions were controlled in three cases which are summarized in \cite{teng2013analytical}.
Five water content probes were injected into the cylinder at depths of $1$, $5$, $10$, $15$ and $19$ $cm$, respectively. The hydraulic parameters of the considered soil are $\alpha=4.8~m^{-1}$, $K_s=3.9\times10^{-6}~m/s$, $\theta_r=0$ and $\theta_s=0.4$. The initial and boundary conditions corresponding to this experiment are:
\begin{equation}
    \begin{cases}
    \Theta(z,0)=1, \\
    \Theta(0,t)=\exp(-\beta t), \\
    \end{cases}
\end{equation}
where $\beta$ is a positive constant given by $\beta=4Db^2$ \cite{teng2013analytical} and $b\leq \alpha/4$ is a fitting parameter. The values of $b$ are $0.413$, $0.489$ and $0.37$ for cases $1$, $2$, and $3$, respectively. The exact solution associated to this numerical test is given by \cite{teng2013analytical}:
\begin{multline}
    \Theta(z,t)=\dfrac{1}{2} \left[ \text{erfc}\left(\dfrac{ D\alpha t-z}{2\sqrt{Dt}}\right)- \e^{\alpha z}\text{erfc}\left(\dfrac{ D\alpha t+z}{2\sqrt{Dt}}\right)\right]+ \\ \dfrac{1}{2}\e^{(\frac{\alpha z}{2}-4b^2Dt)} \left[\e^{2\alpha z} \text{erfc}\left(\dfrac{ z}{2\sqrt{Dt}}+2c\sqrt{Dt}\right)+ \e^{-2\alpha z}\text{erfc}\left(\dfrac{ z}{2\sqrt{Dt}}-2c\sqrt{Dt}\right)\right],
\end{multline}
where erfc is the complementary error function and $c=\sqrt{(\alpha/4)^2-b^2}$.
We display in Figure \ref{evapTest} (left) a comparison of water content distributions at different times of measured, approximate and exact solutions. 

We can observe that the drying rate is higher for the case where the $b$ value is large. It is evident that a higher drying rate results in a lower moisture content profile, whereas a slower drying rate persists for a more extended period. For instance, in the second case where $b=0.489$, the water content takes nearly $200~h$ to reach $10\%$. In contrast, the third case requires a minimum of $300~h$.
The evaporation rates under the three cases are displayed in Figure \ref{evapTest} (right). The computed evaporation rate $E(t)$ by neglecting the gravity is given by \cite{teng2013analytical}:
\begin{equation}
    -E(t)=-D\dfrac{\partial\theta}{\partial z} \Big\vert_{z=0}.
\end{equation}
Scatter points, dashed and solid lines present the experimental profile, the exact and approximate solutions, respectively. A good agreement is observed in comparison to the exact solution. Similarly, a satisfactory correspondence with measured data is obtained, confirming the efficacy of the proposed numerical method in accurately predicting evaporation. For each case, both experimental and numerical results show two distinct evaporation phases. The constant rate stage persists for approximately $90$, $60$, and $180$ hours in cases 1, 2, and 3, respectively. Subsequently, a falling rate phase occurs, followed by the initiation of a residual falling-rate phase at $250$, $200$, and $370$ hours, for cases 1,2 and 3 respectively. The results are obtained using $\varepsilon=0.4$, $n_s=3$, $N_z=200$ and $\Delta t=0.01$. 
The root mean squared errors (RMSE) for water content in the three cases, obtained by comparing the approximate and exact solutions at $T=600~h, 350~h, \text{ and } 900~h$, are $3.8\times 10^{-8}$, $5.31\times 10^{-8}$, and $1.98\times 10^{-8}$ for cases 1, 2, and 3, respectively. These results show the efficacy of the LRBF technique in providing accurate numerical solutions.
The numerical results confirm that the proposed numerical model is able to predict the water content distribution during the evaporation process.

\subsection{Test 2: Unsaturated flow in rooted soils under variable surface flux conditions}
Here, we take into account the presence of the root water uptake in our model, and we use the stepwise \eqref{E7} and exponential \eqref{E8} forms as sink terms in the governing equation. We conduct numerical experiments of the distribution of water content, pressure head, and water flux through a rooted soil with a depth of $L=100~cm$ and a root depth of $40~cm$, which corresponds to $L_1=60~cm$. We consider the following hydraulic parameters: $\theta_s=0.45$, $\theta_r=0.2$, and $K_s=1~cm~h^{-1}$ \cite{SrivastavaYeh,yuan2005analytical}. The rate of reduction in root uptake $\beta$ is set at $0.04$ $m^{-1}$, and the maximum water uptake $R_0$ is set at $0.02~h^{-1}$ for $\alpha=0.01~cm^{-1}$ and $0.0025~h^{-1}$ for $\alpha=0.1~cm^{-1}$, respectively. The following boundary and initial conditions are used:
\begin{equation}
    \begin{cases}
    h(z,0)=h_0(z), \\
    h(0,t)=0, \\
   \left[K(h) \left( \dfrac{\partial h}{\partial z}+1 \right) \right ]_{z=L} =-q_1(t),
    \end{cases}
\end{equation}
where $h_0$ is the initial pressure head and $q_1$ is the time-dependent flux at the soil surface. Both steady-state and transient surface fluxes are considered in this numerical test. 

\begin{figure}[ht!]
\centering
\begin{tabular}{cc}\includegraphics[width=7cm,height=8.4cm,angle=0]{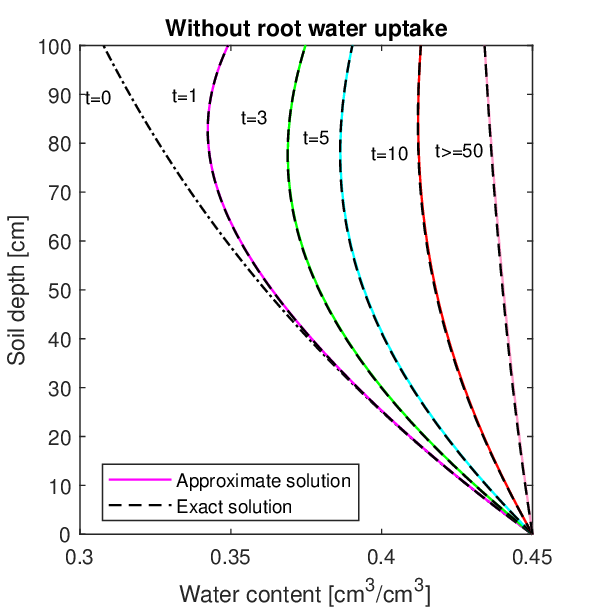}&
\includegraphics[width=7cm,height=8.4cm,angle=0]{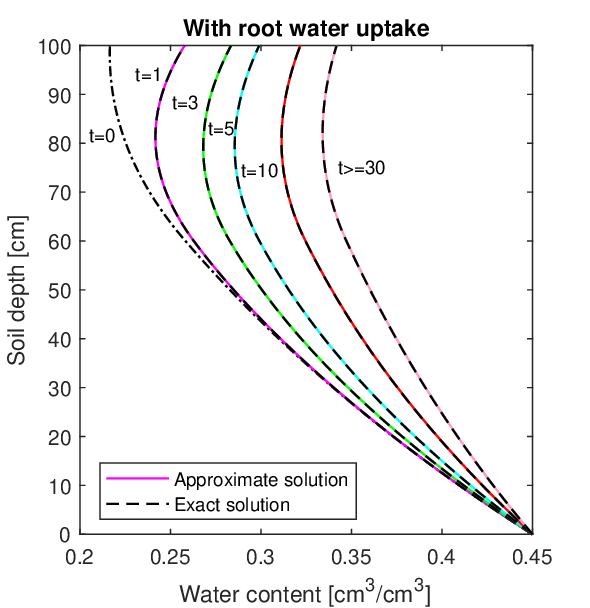}  \\
\includegraphics[width=7cm,height=8.4cm,angle=0]{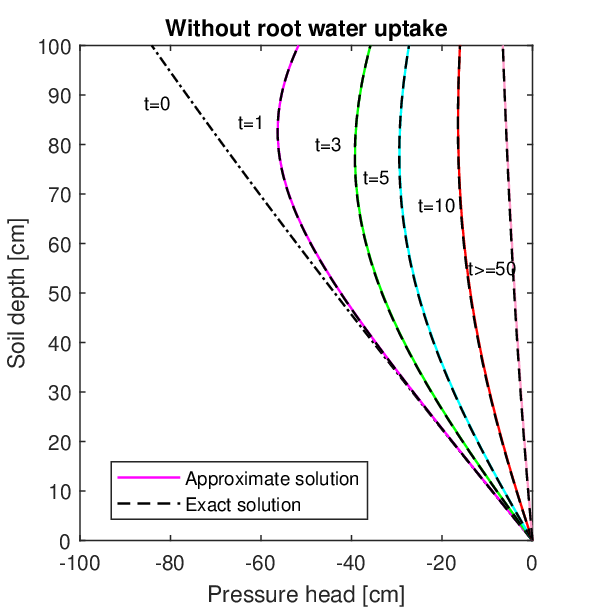} &\includegraphics[width=7cm,height=8.4cm,angle=0]{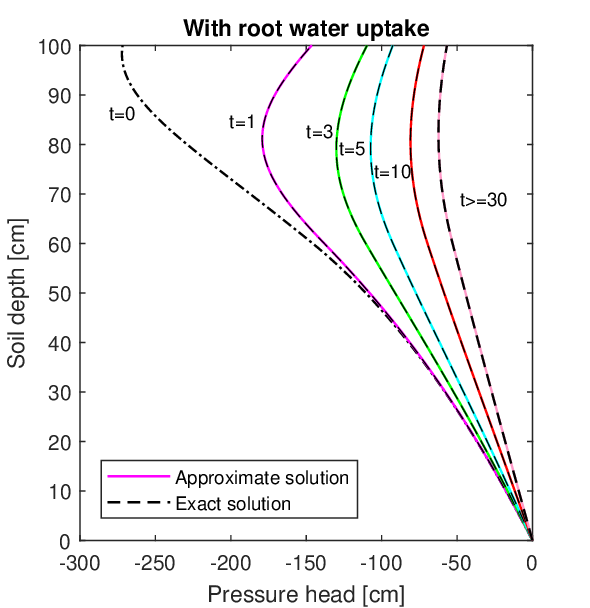}   
\end{tabular}
\caption{Case 1: $\alpha=0.01~cm^{-1}$ and $R_0=0.02~h^{-1}$. Comparison of the water content and pressure head results of the approximate and exact solutions. Left: with root water uptake. Right: without root water uptake.}\label{Hevsite01}
\end{figure}

\begin{figure}[ht!]
\centering
\begin{tabular}{cc}\includegraphics[width=7cm,height=8.4cm,angle=0]{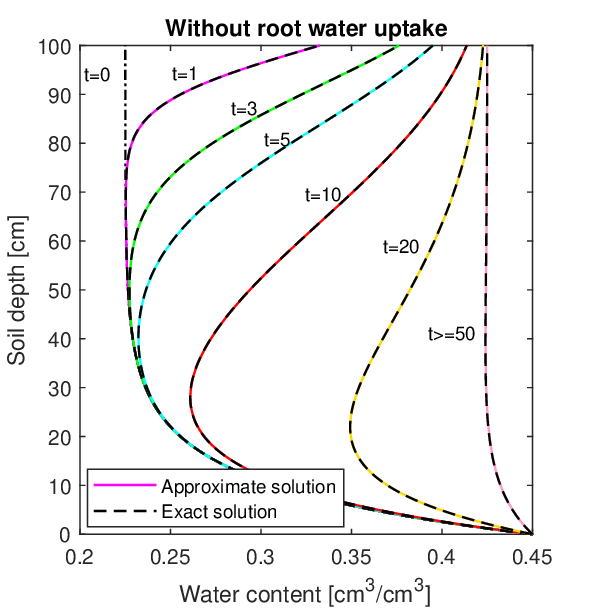}&
\includegraphics[width=7cm,height=8.4cm,angle=0]{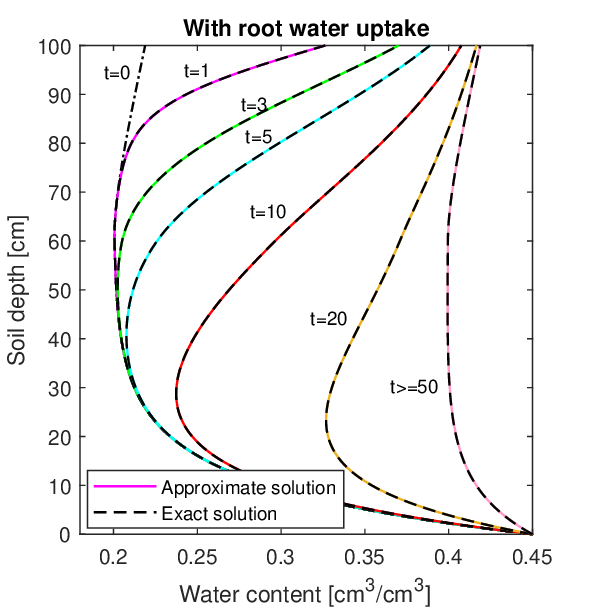}  \\
\includegraphics[width=7cm,height=8.4cm,angle=0]{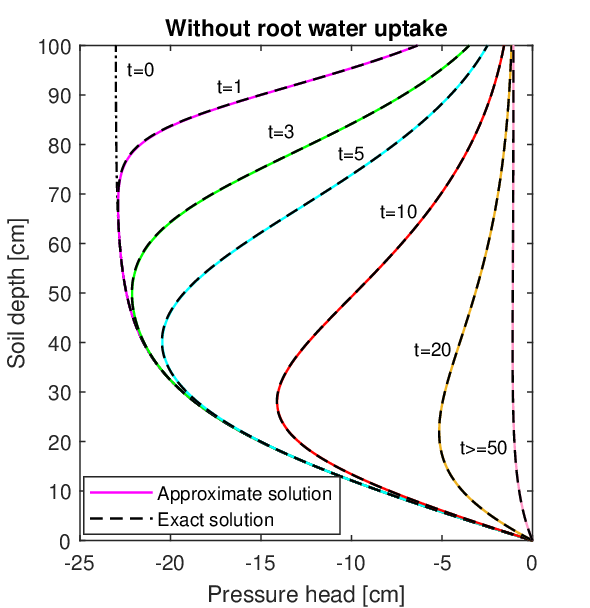} &\includegraphics[width=7cm,height=8.4cm,angle=0]{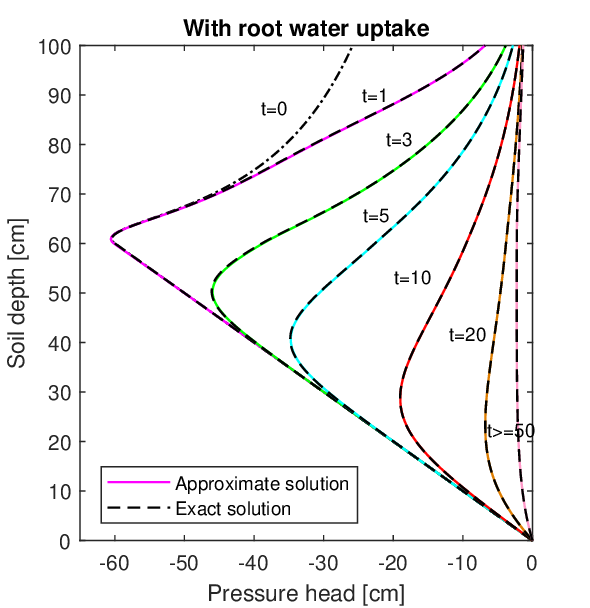}   
\end{tabular}
\caption{Case 2: $\alpha=0.1~cm^{-1}$ and $R_0=0.0025~h^{-1}$. Comparison of the water content and pressure head results of the approximate and exact solutions. Left: with root water uptake. Right: without root water uptake.}\label{Hevsite1}
\end{figure}

We start with steady-state surface flux for the upper boundary condition where we consider a constant infiltration flux of $q_1=-0.9~cm~h^{-1}$. In this case the stepwise \eqref{E5} formulation is used for root water uptake. We display in Figures \ref{Hevsite01} and \ref{Hevsite1} the distribution of the water content and pressure head for $\alpha=0.01~cm^{-1}$ and $\alpha=0.1~cm^{-1}$ respectively. 

The numerical simulations are conducted over a duration of 50 hours. To investigate the impact of water uptake by plan roots on the distribution of soil moisture and pressure head, the time evolution of $h$ and $\theta$ is presented with (right) and without (left) root water uptake.
We can observe that the impact of the root water uptake is significant, particularly when $\alpha=0.01~m^{-1}$ and $R_0=0.02~h^{-1}$. However, the change is relatively weak for the water content in the case of $\alpha=0.1~m^{-1}$ and $R_0=0.0025~h^{-1}$ due to the low value of the maximum water uptake considered. In the first case, the time required to reach a steady state is nearly 50 hours in the absence of root water uptake and 30 hours in the presence of root water uptake. However, in the second case, the time required to attain a stable state is quite similar, at about 50 hours.

\begin{figure}[ht!]
\centering
\begin{tabular}{cc}
\includegraphics[width=8cm,height=6.2cm,angle=0]{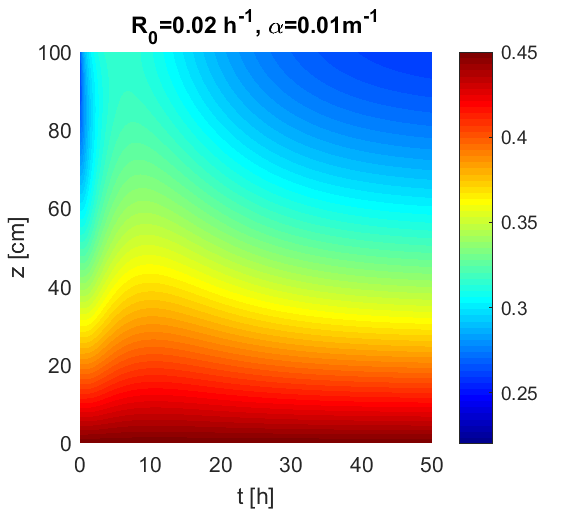} &
\includegraphics[width=8cm,height=6.2cm,angle=0]{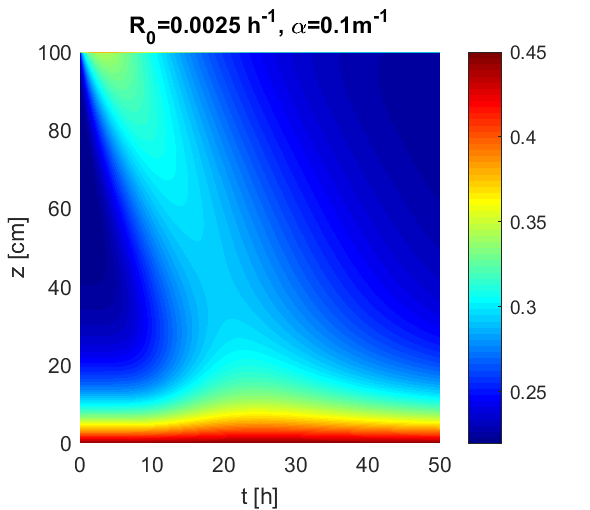}  
\end{tabular}
\caption{Time evolution of soil water content.}\label{ExponT}
\end{figure}

Now, we use a time-dependent surface flux for the upper boundary condition which is more realistic, because the soil surface conditions vary over time due to factors such as evaporation, rainfall and irrigation. 
In this case, the flux at the upper boundary is assumed to decrease exponentially as a function of time, represented by $q_1(t)=q_0+\delta \exp(k_1  t)$, where $\delta=-0.8~cm~h^{-1}$ and $k_1=-0.1~h^{-1}$. The exponential formulation \eqref{E9} is used for root water uptake. The evolution of water content in time and space for the rooted soils is shown in Figure \ref{ExponT}.
The left figure illustrates the rooted soil with $\alpha=0.01~cm^{-1}$ and $R_0=0.02~h^{-1}$. The right figure corresponds to the case where $\alpha=0.1~cm^{-1}$ and $R_0=0.0025~h^{-1}$.
Both soils are supplied with an equal amount of water from the surface, but they exhibit distinct moisture content patterns. The soil profile 1 characterized by $R_0=0.02~h^{-1}$ and $\alpha=0.01~m^{-1}$ is on average wetter than the soil 2 with $R_0=0.0025~h^{-1}$ and $\alpha=0.1~m^{-1}$ even though the volume of water absorbed by the roots is greater than that received from the soil surface.  This difference can be attributed to the soil profile 1 favoring capillary rise, which facilitates water transfer from the water table into the root zones.

\begin{figure}[ht!]
\centering
\begin{tabular}{cc}
\includegraphics[width=7cm,height=5.7cm,angle=0]{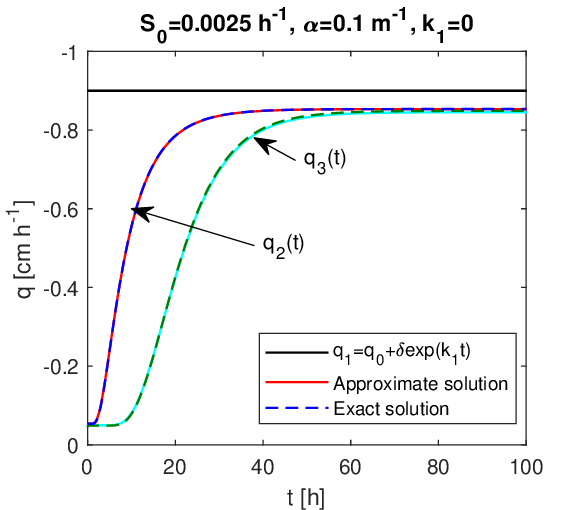} &
\includegraphics[width=7cm,height=5.7cm,angle=0]{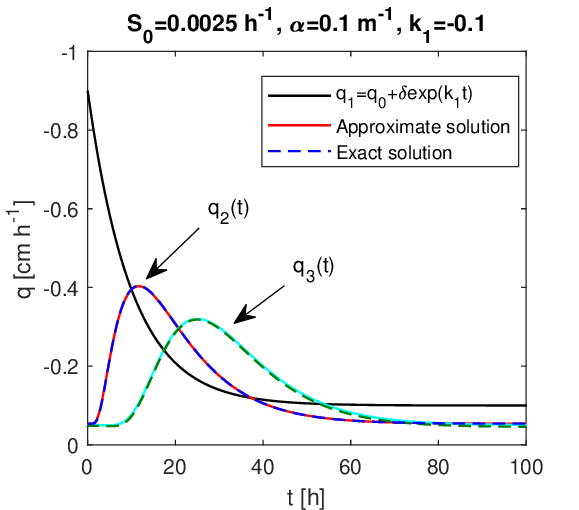} 
\end{tabular}
\caption{Time evolution of water flows at the interface between the root zone and subsoil $q_ 2(t)$ and at the water table $q_ 3(t)$. Left: constant surface flux. Right: varying surface flux.}\label{fluxYuan}
\end{figure}

To assess the transient water flow responding to surface flux changes, we perform numerical simulations of the water flux $q_2(t)$ at the interface between the root zone and subsoil ($z=L_1$), as well as the flow $q_3(t)$ at the water table ($z=0$). Figure \ref{fluxYuan} displays the flows $q_2(t)$ and $q_3(t)$ for both constant (left) and varying surface flux (right).

\begin{figure}[ht!]
\centering
\begin{tabular}{cc}
\includegraphics[width=7cm,height=5.4cm,angle=0]{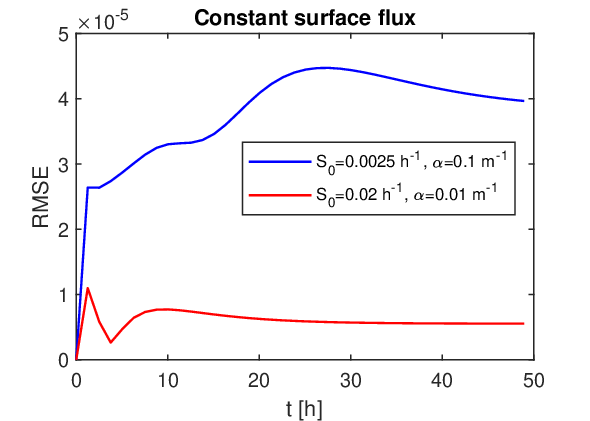} &
\includegraphics[width=7cm,height=5.4cm,angle=0]{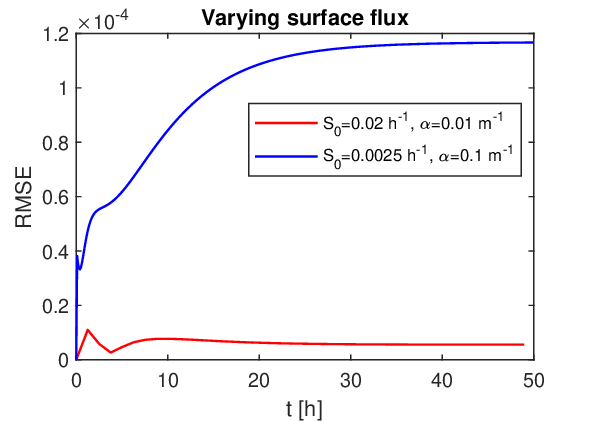} 
\end{tabular}
\caption{RMSE errors of the water content. Left: constant surface flux. Right: varying surface flux.}\label{ErrsYan}
\end{figure}

As the time reaches $50$ hours, $q_2$ and $q_3$ approach $0.5~cm~h^{-1}$ for the constant surface flux and near $0.6~cm~h^{-1}$ for the transient surface flux. The results are obtained using $\varepsilon=0.2$, $n_s=5$, $N_z=1001$ and $\Delta t=0.005$. In all the studied cases, a good agreement is observed between the approximate solutions and the exact solutions which are developed in \cite{yuan2005analytical}. The root mean square errors (RMSE) of the water content for these cases are shown in Figure \ref{ErrsYan}.
The error values are very small for both constant and varying surface fluxes which highlight the accuracy of the proposed LRBF method. This efficiency is particularly important in cases where the surface fluxes may change over time or under different environmental conditions, as it ensures that the proposed method can provide reliable results under a variety of scenarios.
\subsection{Test 3: 2D flow from periodic irrigation furrows}
\begin{figure}[ht!]
\centering
\begin{tabular}{cc}
\includegraphics[width=10cm,height=4cm,angle=0]{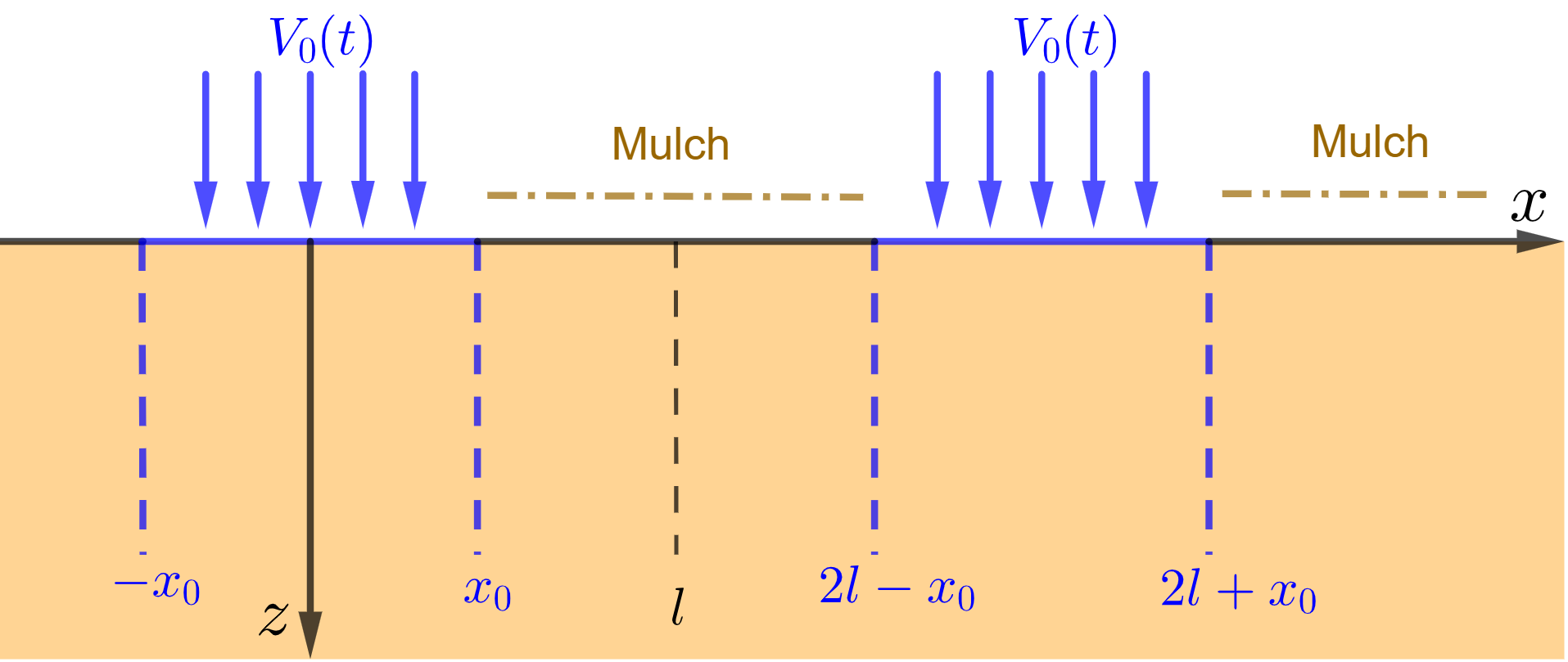}  &
\includegraphics[width=3.9cm,height=3.8cm,angle=0]{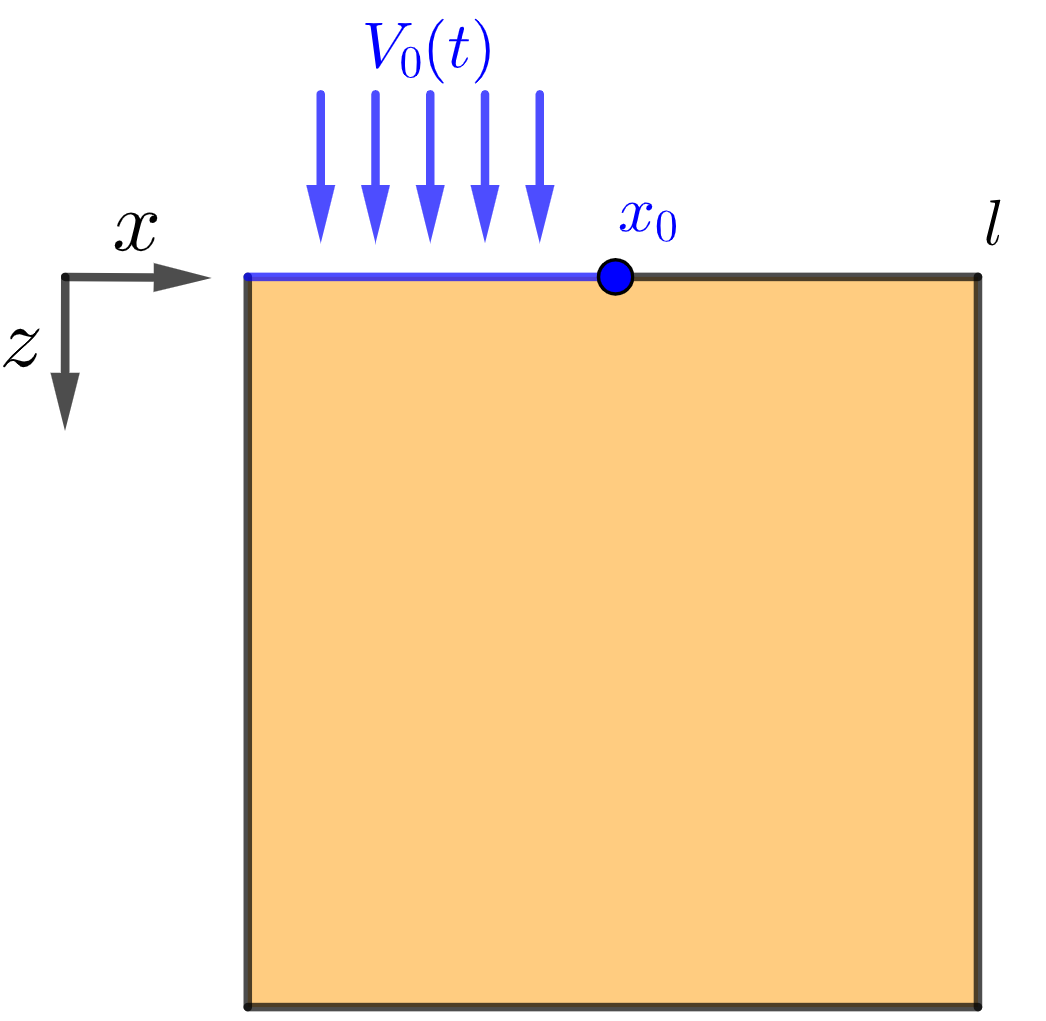}  
 
\end{tabular}
\caption{Schematic of the periodic irrigation furrows (left). The computational domain (right). }\label{f3}
\end{figure}

In this numerical test, we consider surface irrigation furrows with an irrigation rate of $F_0=4~cm/h$. The furrows have a width of $2x_0$ and are periodically placed apart with a spatial period of $2l>2x_0$, as illustrated in Figure \ref{f3} (left). 
The prescribed uniform vertical flux through the furrows $V_0(t)=F_0\e^{(A/t_st)}$ ensures water infiltration into the soil at a specified rate, while the zero flux through the surface between furrows, which is presumed to be protected by mulch, maintains soil moisture by preventing water movement.
These boundary conditions at the soil surface exhibit symmetry due to the periodic placement of the irrigation furrows and the uniform behavior of water infiltration and conservation between these furrows. Owing to the symmetry present in this problem, we can simplify our study by restricting the computational domain to a rectangular region $[0, l]$ as shown in Figure \ref{f3} (right). The horizontal water movement must be absent along the lines of symmetry at $x=0$ and $x=l$, maintaining zero flux across these planes.
\\
We use Model II \eqref{E10} for plant root water uptake, where the sink term depends on soil depth and water content. This dependency increases the nonlinearity of the system, resulting in a more complex numerical test compared to previous tests. We select two distinct plants with their respective parameters \cite{broadbridge2017exact}: Plant 1 with $m=3$, $A=-0.0021$, and $k=-3.31\times 10^{-5}$, and plant 2 with $m=5$, $A=-0.007$, and $k=-3.18\times 10^{-8}$. 

We select the Brindabella silty clay loam soil \cite{perroux1981water,broadbridge1988constant}, characterized by the parameters $\theta_{r}=0$, $\theta_{s}=0.485$, $K_{s}=12~cm~h^{-1}$ and $\alpha=0.142~cm^{-1}$ to simulate water content dynamics, taking into account root water extraction. The chosen values for the parameters $l$ and $x_0$ are $4ls$ and $l_s$ respectively, where $l_s=7~cm$ \cite{broadbridge2017exact}.

\begin{figure}[ht!]
\centering
(a)
\begin{tabular}{ccc}
\includegraphics[width=6cm,height=5.3cm,angle=0]{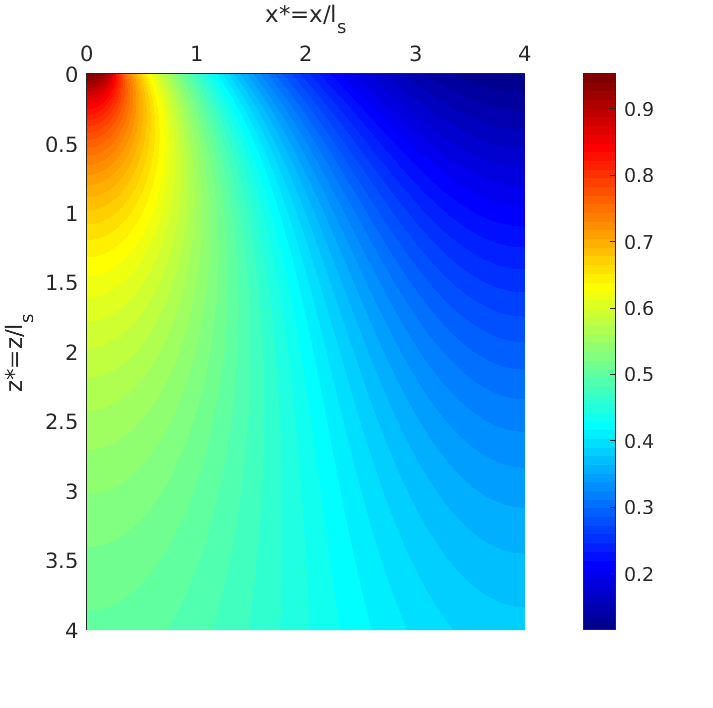}  &
\includegraphics[width=6cm,height=5.3cm,angle=0]{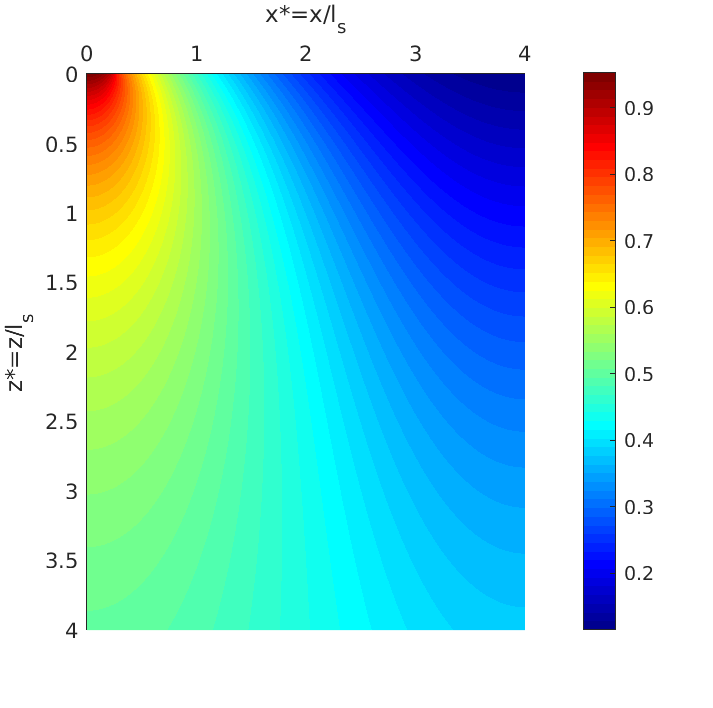} &
\includegraphics[width=5.3cm,height=5.3cm,angle=0]{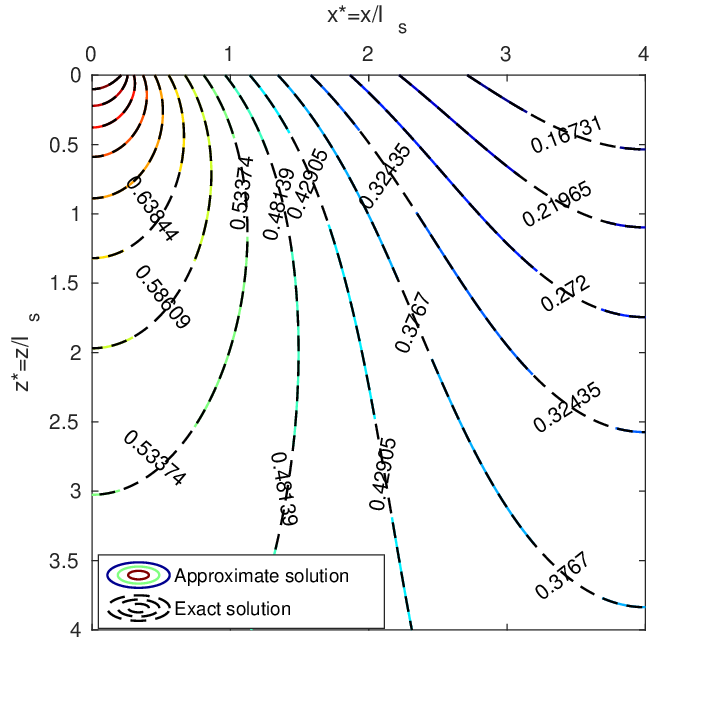}  
\end{tabular}
\centering
(b)
\begin{tabular}{ccc}
\includegraphics[width=6cm,height=5.3cm,angle=0]{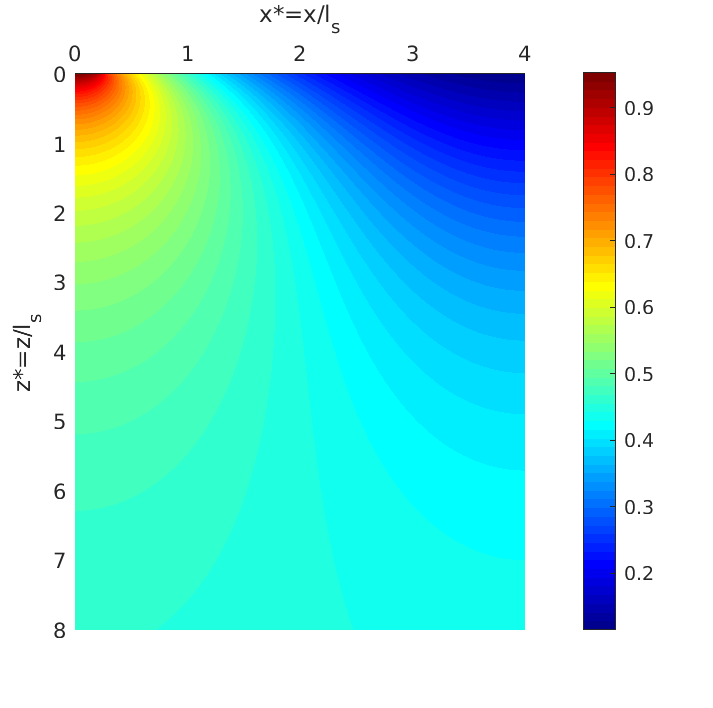} &
\includegraphics[width=6cm,height=5.3cm,angle=0]{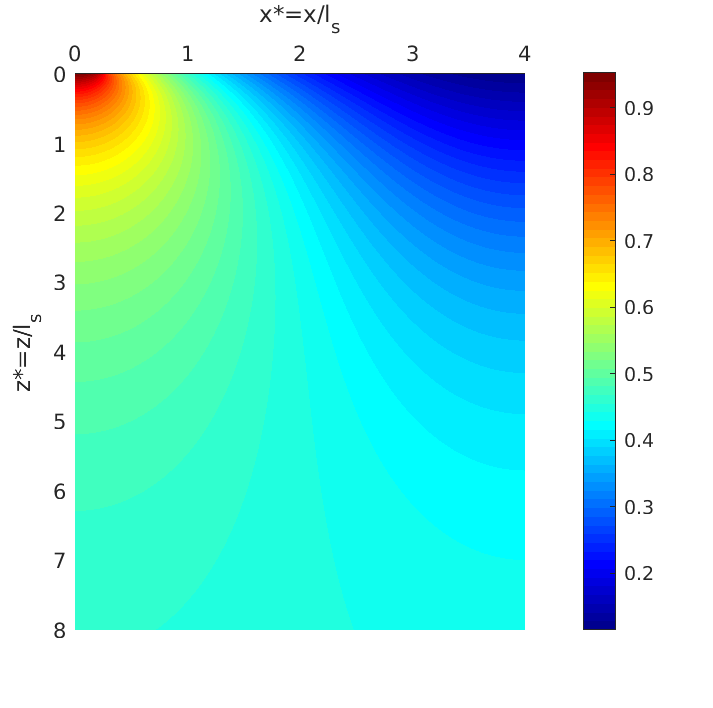}& 
\includegraphics[width=5.3cm,height=5.3cm,angle=0]{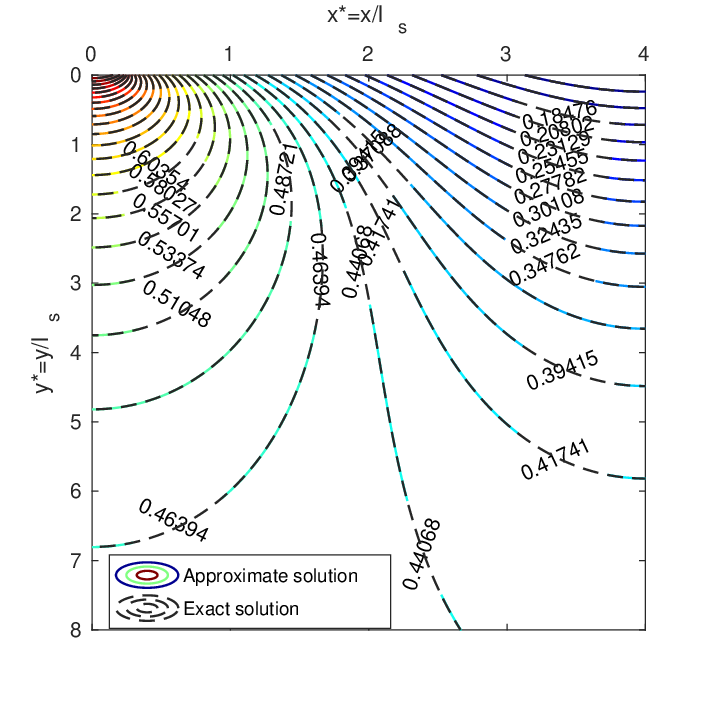} 
\end{tabular}
\caption{Evolution of $\Theta$ near the surface $(z^{*}=4)$ (a), and deep in the soil $(z^{*}=8.5)$ (b): Approximate solution (left), exact solution (middle), and comparison contours (right).}\label{f4}
\end{figure}

The boundary conditions associated to this numerical test can be expressed as follows \cite{broadbridge2017exact}:
\begin{equation}\label{E31}
\begin{cases}
-K\left( \dfrac{\partial h}{\partial z}-1 \right)=V_{0}(t),  & z=0 ~\text{and}~ 0\leq x \leq x_{0},\\
-K\left( \dfrac{\partial h}{\partial z}-1 \right)=0, & z=0 ~\text{and}~ x_{0}< x< l,\\
-K\dfrac{\partial h}{\partial x}=0, & x=0, x=l,\\
-K\left( \dfrac{\partial h}{\partial z}-1 \right)=\e^{(A/t_st)}\left(\dfrac{x_0 F_{0}}{l}\right)\e^{-\alpha(\sqrt{1-4k}-1) z/2}, & z=L.
\end{cases}
\end{equation}
The analytical solution in this case is given by \cite{broadbridge2017exact}:
\begin{equation}
  \mu=\left(\dfrac{\theta_s}{\alpha^2t_s}\right)e^{(A/t_st)} \Phi(x,z),  
\end{equation}
 where, $\Phi$ is given by:
\begin{equation}\label{E32}
\Phi(x,z)=\sum_{j=1}^{\infty}2A_{j}\cos\left( \dfrac{j\pi x}{l}\right)\dfrac{e^{-\alpha(\sqrt{1+4(j\pi l_s/l)^2-4k-1})z/2}}{1+\sqrt{1+4(j\pi l_s/l)^2-4k}},
\end{equation}
where $A_{0}=F_{0} x_{0}/(K_s l)$ and $A_{j}=2F_{0}/(K_sj\pi)\sin(j\pi x_{0}/l)$. According to Equation \eqref{E15}, $\Theta$ can be expressed as follow:
\begin{equation}\label{E33}
\Theta=\dfrac{1}{m}\log\left(1+\alpha^2t_s\left(\dfrac{\e^{m}-1}{\theta_s-\theta_r}\right)\mu\right).
\end{equation}
In our numerical experiments, we set $\varepsilon=0.5$, $\Delta t=0.001$, $n_{s}=5$, $N_{x}=1000$ and $N_{z}=2000$. The evolution of $\Theta$ near the surface $(z^{*}=z/l_s=4)$ and deeper in the soil $(z^{*}=z/l_s=8.5)$, associated with plant 1 and plant 2 respectively, is presented in Figures \ref{f3} and \ref{f4}. These figures illustrate the changes in soil moisture levels at different soil depths as influenced by the distinct root water uptake characteristics of the two plants.

\begin{figure}[ht!]
\centering
(a)
\begin{tabular}{ccc}
\includegraphics[width=6cm,height=5.3cm,angle=0]{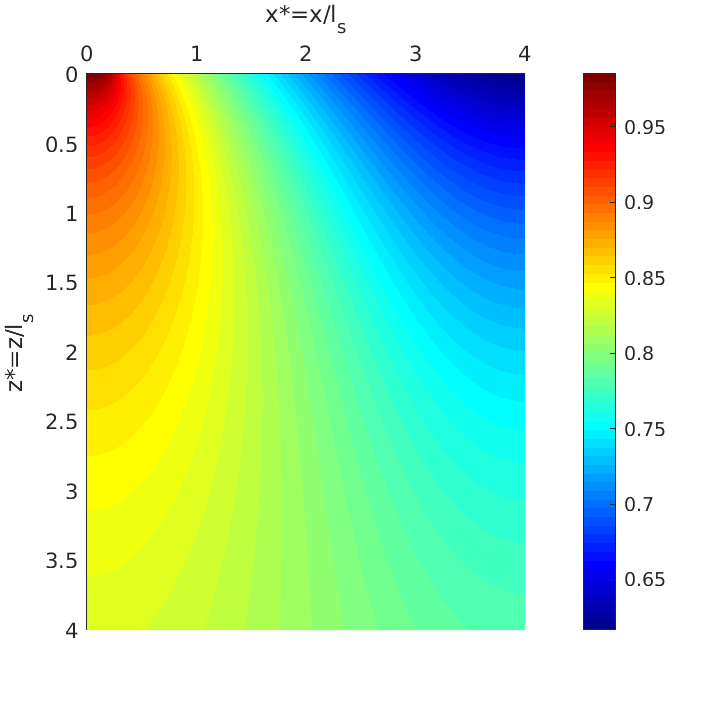} &
\includegraphics[width=6cm,height=5.3cm,angle=0]{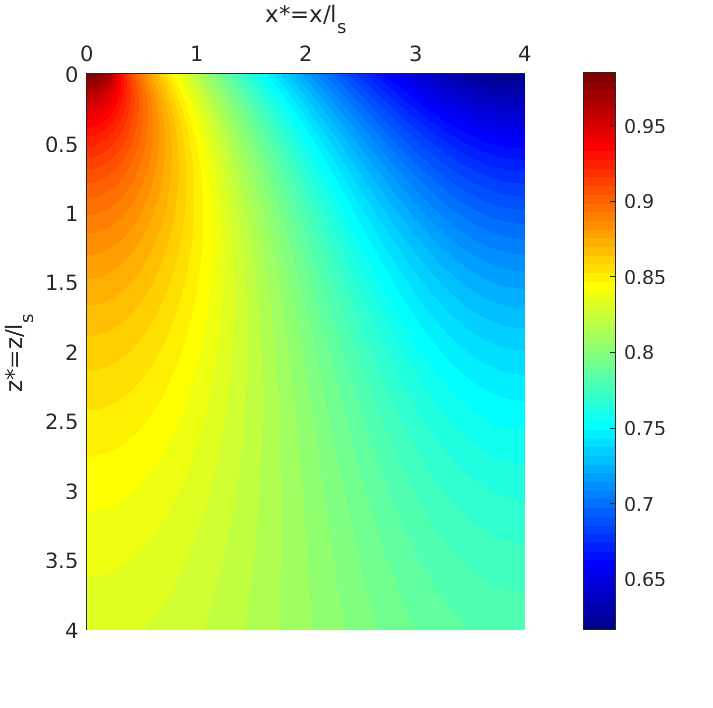} 
\includegraphics[width=5.3cm,height=5.3cm,angle=0]{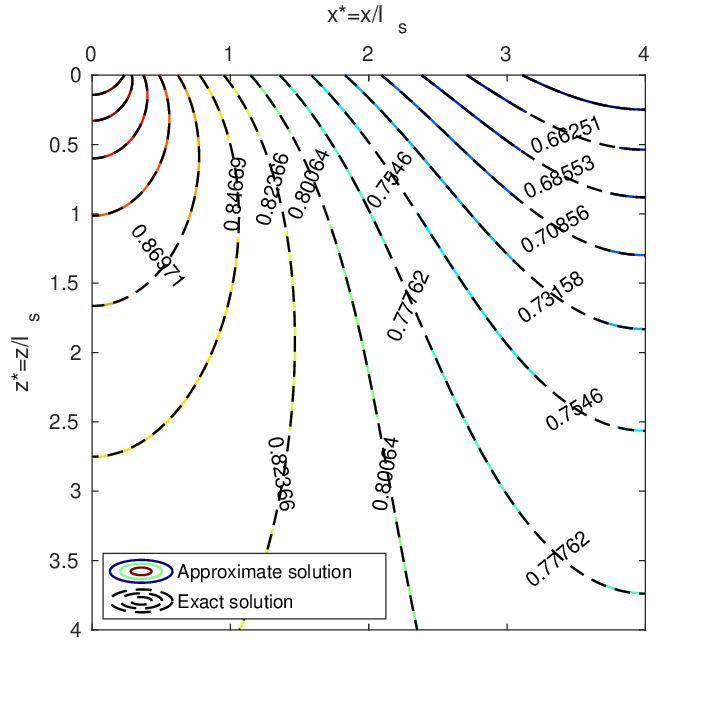} &
\end{tabular}
\centering
(b)
\begin{tabular}{ccc}
\includegraphics[width=6cm,height=5.3cm,angle=0]{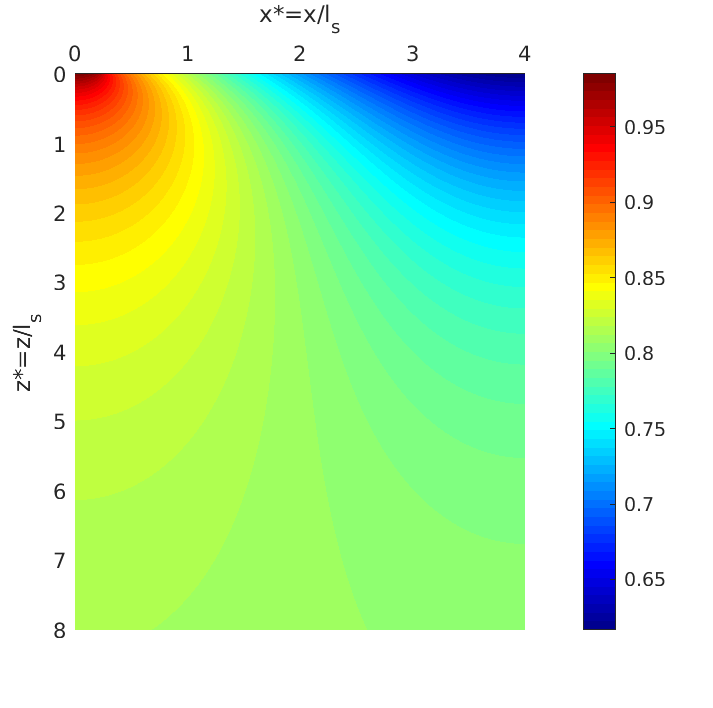} &
\includegraphics[width=6cm,height=5.3cm,angle=0]{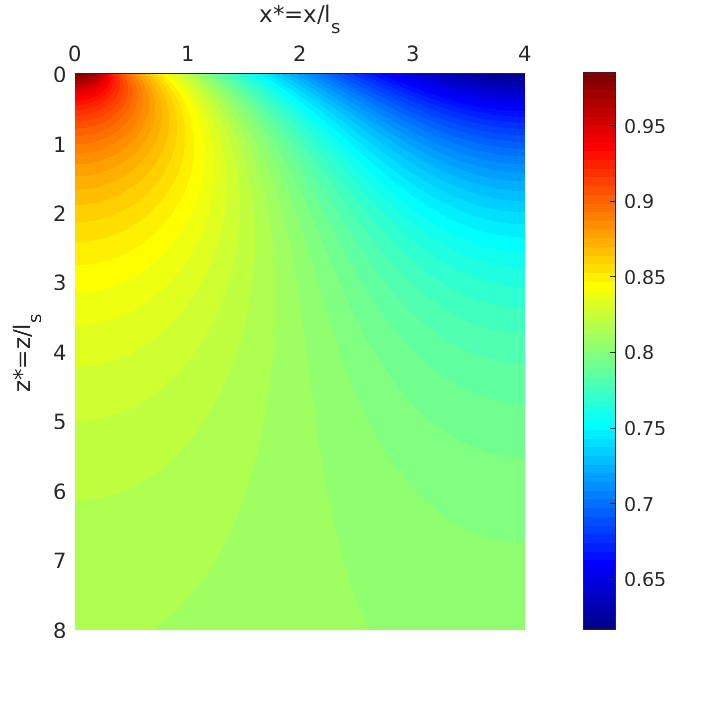} 
\includegraphics[width=5.3cm,height=5.3cm,angle=0]{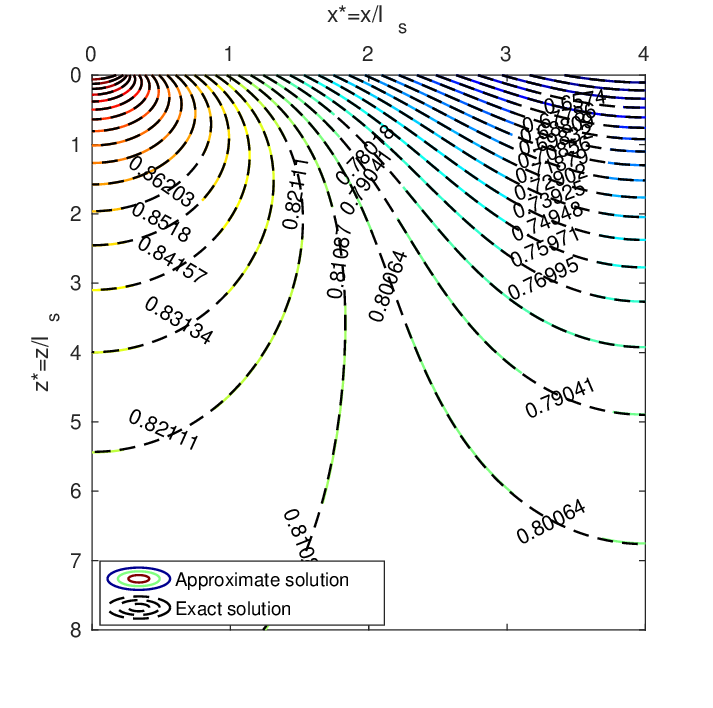} &
\end{tabular}
\caption{Evolution of $\Theta$ near the surface $(z^{*}=1.5)$ (a), and deep in the soil $(z^{*}=9)$ (b): Approximate solution (left), exact solution (middle), and comparison contours (right). }\label{f41}
\end{figure}

We can notice that the soil associated with plant 2 is, on average, more humid than that of plant 1, which shows that the plant corresponding to soil 1 absorbs water at a higher rate than the plant corresponding to soil 2. This observation is confirmed by the values of $\Theta$, which are greater for plant 2 than for plant 1. This difference can be attributed to the distinct values of $k$, $m$, and $A$ chosen for each plant, emphasizing the effect of the plant root system on the soil moisture distribution. For both plant 1 and plant 2, a comparison of the contour plots of $\Theta$ between the approximate and exact solutions is displayed in Figures \ref{f3} and \ref{f4}, respectively. The results show a good correspondence between the two solutions, validating the reliability of the model used in this study.

Table \ref{TE1} presents the $RMSE$ errors between the approximate and exact solutions of $\Theta$ for both plant 1 and plant 2. 
\begin{table}[ht!]
\begin{center}
\caption{RMSE errors between approximate and exact solutions.}\label{TE1}
 \begin{tabular}{|c|c|c|c|}
 \hline
 Plant & $N_{z}$ & $N_{x}$ & $RMSE$ \\
 \hline
 \multirow{2}{1.5cm}{1} & $2000$& $1000$ &$ 6.35\times10^{-5}$   \\
 \cline{2-4}
 & $3000$& $1500$ & $4.23\times10^{-5}$ \\
 \cline{2-4}
 & $4000$& $2000$ & $3.18\times10^{-5}$ \\
 \hline
 \multirow{2}{1.5cm}{2} & $2000$  & $1000$  & $3.73\times10^{-5}$ \\
 \cline{2-4}
 & $3000$  & $1500$ & $2.49\times10^{-5}$  \\

   \cline{2-4}
 & $4000$  & $2000$ & $1.87\times10^{-5}$ \\
  \hline
\end{tabular}
\end{center}
\end{table}
We fix the values of $\Delta t$, $\varepsilon$, and $n_s$, and then vary the number of nodes distributed in the computational domain. The error values, which are small for both plants, diminish as the number of nodes increases. This confirms the convergence and accuracy of the proposed LRBF method in solving the governing equation. The results demonstrate the effectiveness of our approach in simulating soil moisture dynamics while taking into account root water uptake.
 
\subsection{Test 4: 3D irrigation from a circular source}
\begin{figure}[ht!]
\centering
\includegraphics[width=3.6cm,height=6.5cm,angle=0]{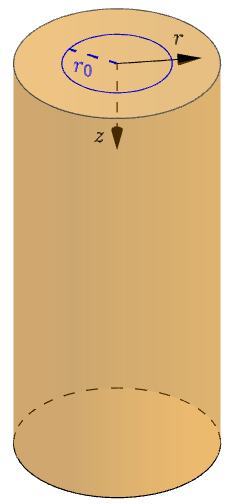}  
\caption{The 3D computational domain. }\label{Fcylinder}
\end{figure}

\begin{figure}[ht!]
\centering
\begin{tabular}{cc}

\includegraphics[width=6cm,height=6.8cm,angle=0]{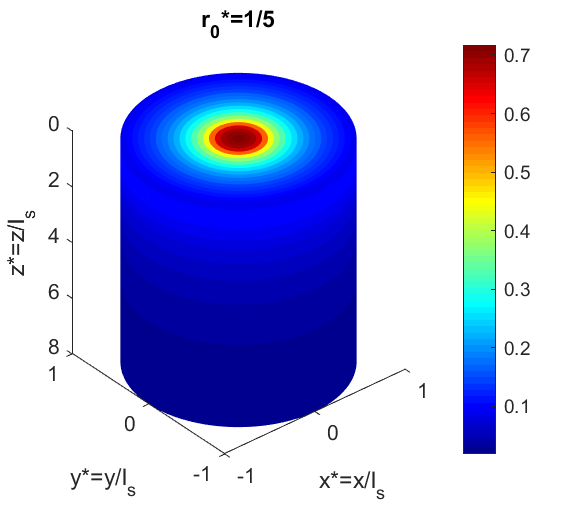}  &
\includegraphics[width=6cm,height=6.8cm,angle=0]{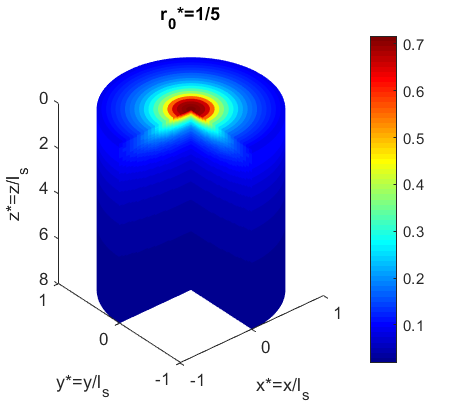} \\
\includegraphics[width=6cm,height=6.8cm,angle=0]{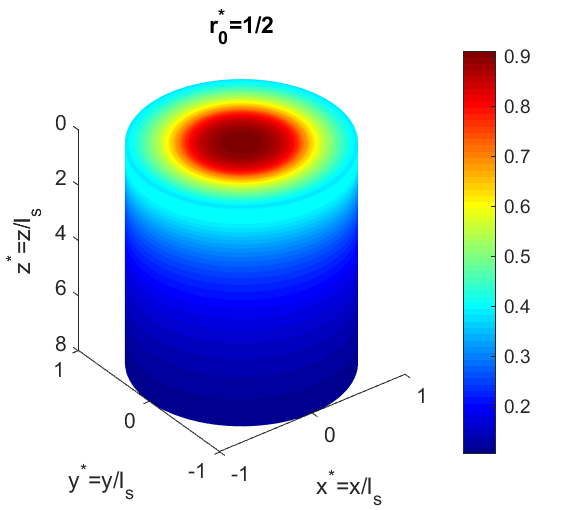}  &
\includegraphics[width=6cm,height=6.8cm,angle=0]{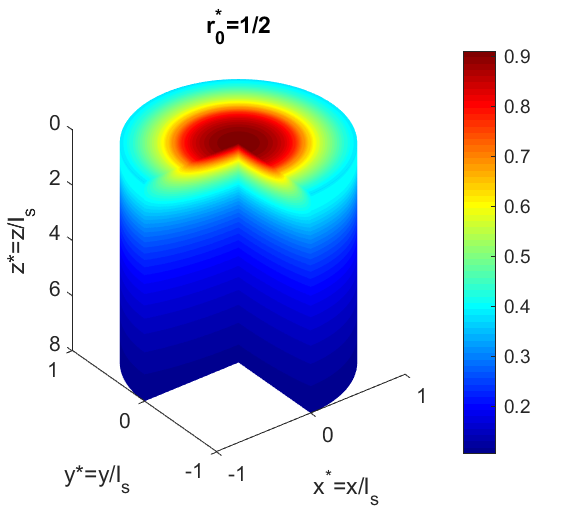} 
\end{tabular}
\caption{Approximate solution of the evolution of $\Theta$ corresponding to plant 1 for $r_0^{*}=1/5$ and $r_0^{*}=1/2$.}\label{F6}
\end{figure}

\begin{figure}[ht!]
\centering
\begin{tabular}{cc}

\includegraphics[width=6cm,height=6.8cm,angle=0]{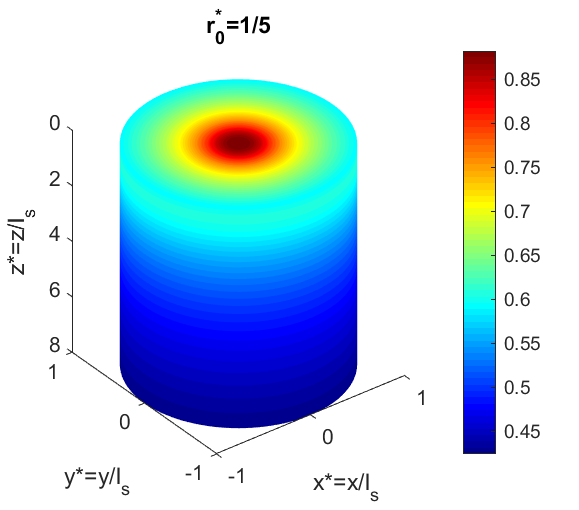}  &
\includegraphics[width=6cm,height=6.8cm,angle=0]{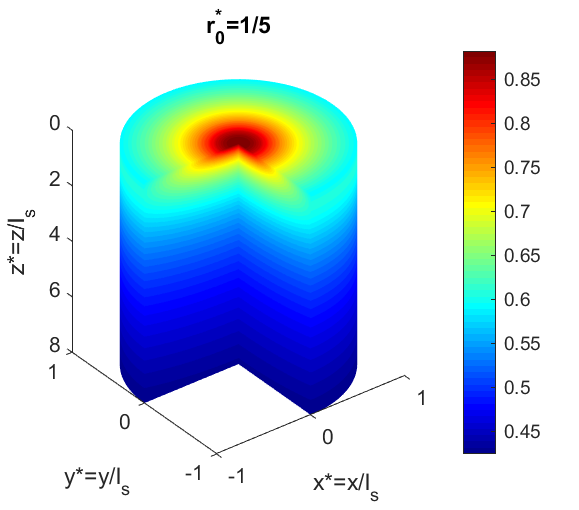} \\
\includegraphics[width=6cm,height=6.8cm,angle=0]{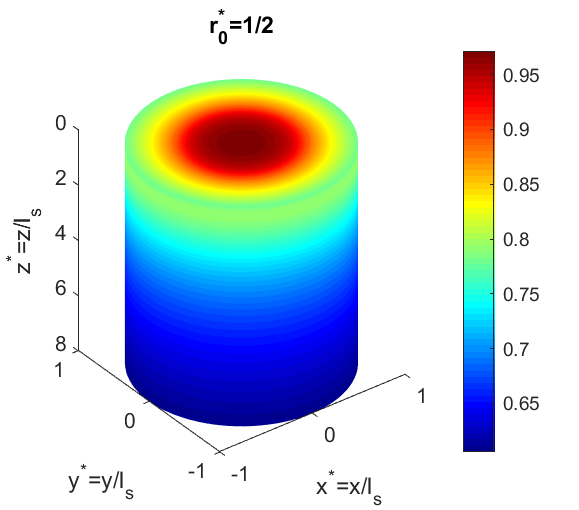}  &
\includegraphics[width=6cm,height=6.8cm,angle=0]{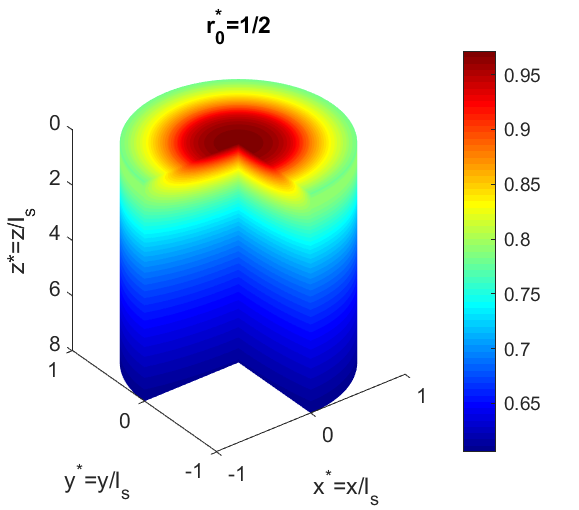} 
\end{tabular}
\caption{Approximate solution of the evolution of $\Theta$ corresponding to plant 2 for $r_0^{*}=1/5$ and $r_0^{*}=1/2$. }\label{F7}
\end{figure}

In this three-dimensional numerical test, we extend the previous analysis on furrow irrigation to assess the effectiveness of the proposed numerical model in predicting 3D soil moisture distribution in the root zone under cylindrical coordinates. We consider a cylindrical region with a length of $L=56~cm$ and a radius of $R=7~cm$ \cite{broadbridge2017exact}, and assume that the irrigation source is a circular region of radius $r_0$ located at the soil surface, as illustrated in Figure \ref{Fcylinder}.
At the top of the cylinder, a prescribed flux boundary $V_0(t)$ is imposed at the circular region of radius $r_0$, and no-flow boundary condition is applied outside of this circular region. Notably, at the lateral sides of the cylinder (i.e., at $r=R$), no-flow boundary condition is imposed. The boundary conditions described above can be written as follows:
\begin{equation}\label{E34}
\begin{cases}
-K\left( \dfrac{\partial h}{\partial z}-1 \right)=V_{0}(t),  & z=0 ~\text{and}~ r  \leq r_{0},\\
-K\left( \dfrac{\partial h}{\partial z}-1 \right)=0, & z=0 ~\text{and}~ r>r_0,\\
-K\dfrac{\partial h}{\partial r}=0, & r=R.
\end{cases}
\end{equation}
We consider model II for plant root water uptake. The analytical solution associated with this numerical test in case of axisymmetric flow is given by \cite{broadbridge2017exact}:
\begin{equation}
  \mu=\left(\dfrac{\theta_s}{\alpha^2t_s}\right)e^{A/t_st} \Phi(r),
\end{equation}
where
\begin{equation}\label{E35}
\Phi=\dfrac{\alpha r_{0}F_{0}}{K_s}\int_{0}^{\infty}\dfrac{2J_{1}(\alpha vr_{0})J_{0}(\alpha vr)}{1+\sqrt{1+4(v^2-k)}}\e^{-\alpha\left(\sqrt{1+4(v^2-k)}-1\right)z/2 }dv,
\end{equation}
where $J_0$ and $J_1$ are Bessel functions of the first kind, for orders $0$ and $1$, respectively. $\Theta$ is calculated using Equation \eqref{E33}.
With the same parameters of soil and plants taken in the previous test, we show in Figures \ref{F6} and \ref{F7} the evolution of saturation for $r_0^{*}=r_0/l_s=1/5$ and $1/2$ for both plants 1 and 2, respectively. 

As in the 2D results, we observe that the soil with plant 2 remains more saturated than that with plant 1, which indicates the significant impact of the corresponding root water uptake parameters on soil moisture distribution. The impact of the radius $r_0$ of the circular irrigation source is also evident from the numerical results. As $r_0$ increases, the soil saturation increases accordingly, which is expected since a larger area is irrigated. This effect is more pronounced near the soil surface, where the irrigation is applied. However, the impact of $r_0$ on the soil moisture distribution decreases with depth as the soil becomes less affected by the irrigation source.

Table \ref{T4} presents the computed errors for $r_0^{*}=1/2$ with varying numbers of collocation points, indicating that the accuracy of the proposed method increases with the increasing of the number of points. This highlights the effectiveness of the numerical model in accurately predicting 3D soil moisture dynamics in the root zone. Moreover, the results demonstrate the capability of the proposed method in handling the complex geometry of the problem and the ability to capture the dynamics of soil moisture distribution under cylindrical coordinates. 
\begin{table}[ht!]
\begin{center}
\caption{RMSE errors between approximate and exact solutions.}\label{T4}
 \begin{tabular}{|c|c|c|c|c|c|c|}
 \hline
 Plant & $N_{z}$ & $N_{r}$ & $\Delta t$ & $\epsilon$ &$n_{s}$ & $RMSE$ \\
 \hline
 \multirow{2}{1.5cm}{1} & $500$& $235$ &  $0.01$ &  $0.5$ & $7$  &$4.94\times10^{-4}$    \\
 \cline{2-7}
 & $1000$& $125$ & $0.01$ &  $0.5$&$7$  & $1.62\times10^{-4}$ \\

 \hline
 \multirow{2}{1.5cm}{2} & $500$  & $62$   &  0.01& $0.3$& $7$& $5.03\times10^{-4}$  \\
 \cline{2-7}
 & $1000$  & $125$   &  $0.01$ & $0.5$& $7$ & $5.83\times10^{-5}$  \\

  \hline
\end{tabular}
\end{center}
\end{table}

The proposed model offers a robust and efficient numerical approach for simulating soil moisture dynamics in root zone, which can be used to inform irrigation management strategies and improve crop yields. The efficiency of the LRBF method makes it a valuable tool for investigating the impact of different irrigation scenarios and soil properties on soil moisture distribution. 
\section{Conclusion}\label{sec:5}
In this paper, we focus on modeling soil moisture distribution in the root zone and introducing computational techniques for efficiently solving the coupled numerical model of infiltration in soils and plant root water uptake. The proposed model is based on the Richards equation and different formulations for root water extraction to study the impact of plant root water uptake on the soil moisture distribution. The Gardner model is used for the capillary pressure. Incorporating an implicit sink term of root water uptake in the system increases the complexity in terms of numerical resolution. The Kirchhoff transformation of the governing equation is employed to reduce the nonlinearity of the system. Our numerical approach is based on LRBF meshless method to solve the coupled numerical model. This technique provides significant computational advantages, including reduced computational cost and accurate numerical solutions since it does not require mesh generation. Our proposed coupled model is validated through various numerical tests, including simulations of soil moisture distribution during evaporation and root water uptake processes in one-, two-, and three-dimensional cases. The reliability of our approach is validated against both experimental data and non-trivial exact solutions. The numerical results show the capability of the proposed method to effectively predict the dynamics of unsaturated flow through soils under evaporation and plant root water absorption. Our study offers a robust numerical framework for studying the impact of root water uptake on soil moisture dynamics, which is important for agriculture and ecosystem management. The proposed numerical techniques could be further extended to incorporate more complex root water uptake models, which can provide a promising avenue for future research.
\section*{\textbf{Acknowledgments}}
AB gratefully acknowledges funding from UM6P-OCP. AT and AB gratefully acknowledge funding for the APRD research program from the Moroccan Ministry of Higher Education, Scientific Research and Innovation and the OCP Foundation.

\bibliographystyle{elsarticle-num}
\bibliography{mabiblio.bib}

\end{document}